\documentclass[11pt]{amsart}
\usepackage{amsmath}
\usepackage{amsfonts}
\usepackage{amssymb}
\usepackage{amscd}
\usepackage{epsfig}

\catcode `\@=11
\def\numberbysection{\@addtoreset{equation}{section}
         \renewcommand{\theequation}{\thesection.\arabic{equation}}}
\numberbysection
\def\subsubsection{\@startsection{subsubsection}{3}%
  \normalparindent{.5\linespacing\@plus.7\linespacing}{-.5em}%
  {\normalfont\bfseries}}
\newtheorem{thm}{Theorem}[section]
\newtheorem{lem}[thm]{Lemma}
\newtheorem{prop}[thm]{Proposition}

\newtheorem{introthm}{Theorem}
\newtheorem{introprop}{Proposition}

\theoremstyle{definition}

\newtheorem{df}[thm]{Definition}

\newtheorem{rmk}[thm]{Remark}
\newtheorem{nota}[thm]{Notation}

\def\nn{\nonumber}

\def\In{In}

\def\M{\mathcal{M}}

\def\Sn{\mathbb{S}_n}
\def\Sm{\mathbb{S}_m}

\def\Z{\mathbb{Z}}

\def\del{\partial}

\def\t{\tau}
\def\G{\Gamma}
\def\a{\alpha}

\def\Cact{\mathcal{C}act}

\def\CWcact{K}

\def\stinf{\pi_{\infty}}
\def\stinftop{\stinf^{\text{\it top}}}
\def\ininf{i_{\infty}}

\newcommand\val[1]{\text{\it val}(#1)}
\newcommand{\arity}[1]{|#1|}

\newcommand{\app}[1]
{\addtocounter{section}{1}
\section*{Appendix \thesection : #1}
}

\def\color{{clr}}

\def\lab{{Lab}}

\def\DArc{\mathcal{DA}rc}

\def\Arc{\mathcal{A}rc}

\def\Ainf{A_{\infty}}
\def\Cinf{\Cact_{\infty}}

\def\stable{{\mathcal T}_{\infty}}
\def\bipartite{{\mathcal T}_{\text{\it bipart}}}
\def\heights{{\mathcal T}_{\text{\it ht}}}
\def\topheights{\heights^{\text{\it top}}}
\def\plantree{{\mathcal T}_{\text{\it pp}}}
\def\corollas{{\mathcal T}_{\text{\it cor}}}
\def\Tass{\plantree}
\def\Tassheight{\Tass^{\text{\it ht}}}
\def\Tcyc{{\mathcal T}_{\text{\it cyclo}}}
\def\Tcycheight{\Tcyc^{\text{\it ht}}}

\def\Rp{{\mathbb R}_{>0}}
\def\R{{\mathbb R}}
\def\height{{h}}
\def\depth{depth}
\def\topheight{{w}}
\def\Cyc{{\text{\it Cyc}}}
\def\Cyctop{\Cyc^{\text{\it top}}}
\def\Linf{\mathcal{L}in\mathcal{T}ree_{\infty}}

\def\var{v}
\def\varedges{E_{\var}}
\def\Evar{\varedges}

\def\codim{codim}

\def\CWheight{K^{\text{\it ht}}}
\def\CWass{K^{\infty}}
\def\CWcact{K^1}
\def\CWcyc{K^{\text{\it cyc}}}
\def\retract{r}
\def\retractchain{r_*}

\def\CCass{CC_*(\CWass)}
\def\CCcact{CC_*(\CWcact)}

\def\wangles{\angle^{w}}

\def\Emix{E_{\text{\it mixed}}}
\def\Eblack{E_{\text{\it black}}}
\def\Ewhite{E_{\text{\it  white}}}
\def\Vblack{V_{\text{\it  black}}}
\def\Vwhite{V_{\text{\it white}}}

\def\starwhite{*^w}

\def\cell{C}
\def\E{{\mathcal E}}

\def\intop{\imath}

\def\In{In}
\def\Out{Out}

\def\Dsul{{\mathcal DS}ul}
\def\CWsul{\Dsul^1}

\begin{document}

\title[Associahedra, Cyclohedra and the
$A_{\infty}$--Deligne
conjecture]
{Associahedra, Cyclohedra and a  Topological
 solution to the $A_{\infty}$--Deligne
conjecture}

\author
[Ralph M.\ Kaufmann]{Ralph M.\ Kaufmann}
\email{rkaufman@math.purdue.edu}

\address{Purdue University, Department of Mathematics, 150 N. University St.,
West Lafayette, IN 47907--2067}

\author
[Rachel Schwell]{Rachel Schwell}

\address{Trinity College,
Department of Mathematics
300 Summit St.,
Hartford, CT 06106--3100
}
\email{Rachel.Schwell@trincoll.edu}

\begin{abstract}
We give a topological solution to the $\Ainf$ Deligne conjecture
using   associahedra and cyclohedra. For this we construct  three CW
complexes whose cells are indexed by products of polytopes.
Giving new explicit realizations of the polytopes in terms of different
types of trees, we are able to show that the CW complexes are cell models
for the little discs. The cellular chains of one
complex in particular, which is built
out of associahedra and cyclohedra, naturally acts on the Hochschild cochains
of an $\Ainf$ algebra  yielding an explicit,  topological and minimal
solution to the  $\Ainf$ Deligne conjecture.

Along the way we obtain new results about the cyclohedra, such
as a new decompositions into products of cubes and simplices,
which can be used to
realize them via a new iterated blow--up construction.
\end{abstract}

\maketitle


\section*{Introduction}
In the last years Deligne's conjecture has been a continued source
of inspiration. The original conjecture  states that there is a
chain model of the little discs operad that acts on the Hochschild cochains
of an associative algebra, which induces the known Gerstenhaber structure
\cite{G} on cohomology. It has by now found many proofs,
\cite{Maxim,T,MS,Vor2,KS,MS2,BF,del}
which all have their unique flavor.
This plethora of approaches comes from
the freedom of choice of the chain model for the little discs operad.
Among these there are
``minimal'' choices which are cellular and have exactly the cells one needs
to give the relevant operations induced by the operadic structure
\cite{MS,MS2,del}.
In the $\Ainf$ algebra setting where one only assumes that the
algebra is homotopy associative, astonishingly there has so far been only one
solution \cite{KS} based on homological algebra,
  although this subject if of high current
 interest for instance in  Mirror--Symmetry, the theory of $D$--branes
and String Topology.

 In this paper,
we give a new topological, explicit, ``minimal'' solution via
a cell model
for the chains of the little discs which acts on the Hochschild complex
of an $\Ainf$ algebra. This is the geometrization of the combinatorial
Minimal Operad $\M$ introduced by Kontsevich and Soibelman \cite{KS}.

\begin{introthm}[Main Theorem]
\label{theorem1}
There is a cell model $\CWass$ for the little discs operad, whose operad
of cellular chains acts on
the Hochschild cochains of an $\Ainf$ algebra inducing
the standard operations of the homology of the little discs operad
on the Hochschild cohomology of the algebra.
Moreover this cell model is minimal in the sense
that  the cells correspond exactly to the natural operations obtained
by concatenating functions and using the $\Ainf$ structure maps.
\end{introthm}

This statement is a statement over $\Z$.
The first observation which leads us to the proof
 is that the differential of $\M$
is captured by the combinatorics of associahedra and cyclohedra.
This allows us to construct a CW model $\CWass$ whose cellular chains are
naturally isomorphic to $\M$.
The proof that this cellular chain operad is a model of chains
for the little discs operad is a bit involved.
For this we need to compare three CW complexes, each of them built
on  polytopes.  The first, $\CWcact$,  is
the cell model of the little discs
which is the on given by normalized spineless cacti \cite{del},
here the polytopes are just simplices. The second
is the cell model $\CWass$ mentioned above; the cells in this complex are
products of associahedra and cyclohedra. And lastly $\CWheight$
which is a mediating cell model constructed from trees with heights.
In this  model the cells are products of cubes and simplices.
There is a chain of five propositions which leads to the Main Theorem:

\begin{introprop}
\label{prop1}
As chain operads $\CCass$ and $\M$ are equivalent.
\end{introprop}

\begin{introprop}
\label{prop2}
The cell models $\CWheight$ and $\CWass$ have the same realization.
Moreover, $\CWheight$ is just a cellular subdivision of $\CWass$.
\end{introprop}

\begin{introprop}
\label{prop3}
The space $|\CWcact|$ is a strong deformation
 retract of $|\CWheight|$.
\end{introprop}

\begin{introprop}
\label{prop4}
The map induced by the retract $\retract:|\CWass|=|\CWheight|\to |\CWcact|$
on the chain level, $\retractchain:\CCass\to \CCcact$, is a
morphism operads. In fact, it is the map $\stinf$ of \cite{del}.
\end{introprop}

\begin{introprop}\cite{KS}
\label{prop5}
$\M$ acts on the Hochschild complex of an $\Ainf$ algebra in the appropriate
fashion, that is it induces the Gerstenhaber structure on the
Hochschild cohomology.
\end{introprop}
The fact that $\M$ acts is true almost by definition; this is presumably why
it is called the ``minimal operad'' in \cite{KS}.

\begin{proof}[Proof of the Main Theorem]
By Proposition \ref{prop1} and \ref{prop5}
we see that  $CC_*(\CWass)$ acts in the appropriate fashion.
By Propositions \ref{prop2} and \ref{prop3} $|\CWass|$
is homotopy equivalent to $|\CWcact|$ and since by \cite{del}  $\CWcact$
is  a CW model
for the little discs, so is $\CWass$. {\em A priori} this only has
to be true on the space/topological level, but by
Proposition \ref{prop4} on homology the retraction map
$\retract$ is an operadic isomorphism and hence $\CWass$
is an operadic cell model.
\end{proof}

This actually
answers a question of Kontsevich--Soibelman \cite{KS} about a
smooth cell model for $\M$. In terms of a CW complex which is minimal in the
above sense it cannot be had. There is however a certain thickening of cells,
which indeed is a smooth manifold model \cite{manfest}.
This is again given by a CW complex defined by trees, but
with slightly different combinatorics.
In this manifold model,  the action on Hochschild is, however, not minimal;
its dimension is already too big.
It is nonetheless a very natural geometric realization and
nicely linked to the arc complex and the Arc Operad of \cite{KLP}.

Our main tool for constructing the CW complexes
are trees. In each case, we fix
a particular combinatorial class of trees with a differential on the free Abelian
group they generate. Based on
this combinatorial data we build CW complexes, which are indexed by the
particular type of tree such that the tree differential gives the gluing maps and
hence we obtain an isomorphism of Abelian groups between the cellular chains
and the Abelian group of trees. The individual cells are assembled out of products of polytopes. These vary depending on the CW model we are constructing as mentioned above. The building blocks we use for $\CWcact,\CWass$ and $\CWheight$
respectively are simplices,
associahedra and cyclohedra, and simplices and cubes.
The operad structures we consider are all induced from the topological level.
In all three cases, pushing the operad structure
to the homology yields an operad isomorphic to the homology
of the little discs operad.

\begin{introthm}
\label{theorem2}
The  realizations $|\CWass|\simeq |\CWheight|$ and $|\CWcact|$ are
all topological quasi--operads and sub--quasi--PROPs of
the Sullivan--quasi--PROP $\CWsul$ of \cite{hoch1}.
There is also a renormalized quasi--operad structure such that the
induced quasi--operad structure on their cellular chains
$CC_*(\CWass)\simeq\Z\stable, CC_*(\CWheight)\simeq\Z\heights$
and $CC_*(\CWcact)\simeq\Z\bipartite$ is an operad structure
and coincides with the combinatorial operad structure on the trees.
Moreover, all these operad structures are models for the little discs operad.
\end{introthm}

The reader familiar with
these constructions of \cite{KLP} and \cite{hoch1}
may appreciate that the gluings here are just tweakings of
the usual gluings of foliations.
In fact, as far as these structures are concerned
the language of arcs on surfaces would be much easier.
 In the main text we phrase everything in
the equivalent language of trees in lieu of that of arcs
since it is a more widely spoken language
and the tree description is needed to define operations
on the Hochschild complex.
We will however provide a short dictionary in Appendix \ref{arcapp} and
 relegate the proof of Theorem \ref{theorem2} and Proposition
\ref{prop1} to this appendix as they are not absolutely essential
to the argument of the Main Theorem. Proposition \ref{prop1} can
be replaced by the {\it ad hoc} Definition \ref{operadassdef} (see
Proposition \ref{textpropone}).

Appendix \ref{arcapp} will be key in providing the $\Ainf$
generalization of the results of \cite{hoch1,hoch2} and hopefully
shed light on the different constructions stemming from
string topology and mirror symmetry providing similar actions. We would like
to emphasize that in the present
study the CW complexes provided by arcs {\em
do not give the chain model that acts directly}
in contrast to the previous
constructions \cite{del,cyclic,hoch2}
where the arc picture directly gave cells that could be used
for the action. Now, for the first time,
we  need
to consolidate the cells into bigger super--cells in order to have
an action, as the original cells are too fine. This realization and
the presented construction are hence essential to the further study
of chain level actions. One other particularly interesting
issue is the renormalization of the quasi--PROP composition.
This is a novel feature that is necessary to obtain the correct
combinatorics for the $\Ainf$ case on the cell level. These  cannot be handled
by the arguments of \cite{hoch1} alone.

In the process of comparing the models, we establish new facts about
the classical polytopes such as the cyclohedron, which are
interesting in their own right.

\begin{introthm}
\label{theorem3}
There is a new decomposition of the cyclohedron $W_{n+1}$ into a simplex and cubes.
Correspondingly, there is an iterated ``blow--up'' of the simplex to a
cyclohedron, with $n-1$ steps. At
 each stage $k$ the polytopes that are glued on are a
 product of a simplex $\Delta^{n-k}$ and a cube $I^k$,
where the factors $\Delta^{n-k}$ attach to
the codimension $k$--faces of the original simplex.
\end{introthm}

So as not to perturb the flow of the main text, Theorem \ref{theorem3}
and details about the cyclohedron that are
not needed in the proof of the Main Theorem
are referred to Appendix \ref{cycloapp}.

The organization of the paper is as follows:

We start by giving the combinatorial background and introducing the relevant types
of trees in \S\ref{treesection}. Here we also discuss
the three operads of Abelian groups with differentials on which the CW
models are based.
Before introducing said models, we turn to the polytopes that will be used
to construct them: simplices, associahedra and cyclohedra
in \S\ref{polysection}. Here we give two CW decompositions
each of the associahedron
and the cyclohedron. The second CW composition is novel and
leads to Theorem \ref{theorem2}.
Armed with these results we construct the three relevant CW complexes
in \S\ref{cellsection} and prove their relations as expressed in Propositions
\ref{prop2}--\ref{prop4}; these are Propositions
 \ref{refinementprop}, \ref{retractprop} and \ref{textprop4}.
In the final paragraph of the main text, \S\ref{ainfdelsection},
we assemble the results to prove the Main Theorem, Theorem \ref{cycthm}.

 Appendix \ref{arcapp}  gives the relationship
to the arc operad and the Sullivan quasi--PROP, and provides
the proofs of Theorem \ref{theorem2} (Theorem \ref{texttheorem2})
and Proposition \ref{prop1}
which using Definition \ref{operadassdef} is Proposition \ref{textpropone}.
Finally, in Appendix \ref{cycloapp}, we distill the results
on the cyclohedron of
the main text to give the sequential blow--up of Theorem
\ref{theorem3} (Theorem \ref{texttheorem3}) and demonstrate
this on the examples of $W_3$ and $W_4$.

\section*{Acknowledgments}
It is a pleasure to thank J.~Stasheff, S.~Devadoss
and J.~McClure for interesting
and useful discussions. R.K.\ would also like to thank the
Max--Planck--Institute for Mathematics in Bonn, Germany, for
its kind hospitality and support.

\section{Trees, dg--Operads and Algebras}
\label{treesection}
\subsection{Trees}
Let us first recall the standard definitions and then fix
the specific technical conditions on the trees
with which we will be working.

A {\em graph} will be a 1-dim CW complex and a {\em tree} will be a
graph whose realization is contractible. We will need some further
data. To fix this data, we note that given a graph $\G$ the set of
$0$--cells forms the set of {\em vertices} $V(\G)$ and the set of one
cells form the set of {\em edges} $E(\G)$. A {\em flag} is a half
edge. The set of all flags is denoted by $F(\G)$. Notice that there
is a fixed point free involution $\imath:F(\G)\mapsto F(\G)$ which
maps each half--edge to the other half edge making up the full edge.
Each flag has a unique vertex, which we will call the vertex of the
flag. The respective map taking a flag to its vertex will be called $\del$.
The flags at a vertex $v$ are the half edges incident to that
vertex. The set of these flags will be denoted by $F_v(\G)$.
The {\em valence} of a vertex $v$ is
defined to be $\val{v}=|F_v(\G)|$.

For us a {\em ribbon graph} is a graph $\G$ together with a cyclic
order on each of the sets $F_v(\G)$. We impose  no condition
on the valence of a vertex. The cyclic orders give rise to a map $N$
which assigns to a flag $f$ the flag following $\imath(f)$ in the
cyclic order. The iteration of this map gives an action of $\Z$ on the set of flags. The
{\em cycles} are the orbits of this latter map.

An {\em angle} $\a$ of a ribbon graph is a pair of flags $\{f_1,f_2\}$
 which share the same vertex $\del(f_1)=\del(f_2)$ and
where $f_2$ is the immediate successor of $f_1$.
Notice that these may coincide.
The edges of $\a$ are $e_i=\{f_i,\imath(f_i)\}$.
There is a 1--1 correspondence between flags (or edges)
at a vertex and the angles
at a vertex.

A ribbon graph is called {\em planar} if its image can be embedded in the plane
in such a way that the induced cyclic orders
coming from the orientation of the plane
equals the given cyclic order of the graph.

A globally marked ribbon graph is a ribbon graph with a
distinguished flag. A globally marked planar tree is traditionally
called {\em planar planted}. In the tree case, the vertex of the marked flag is
called the {\em  root} and denoted by $v_{\rm root}$; the vertices $v$ with
$\val{v}=1$ which are distinct from $v_{\rm root}$
will be called {\em leaves} and the set of these vertices will be denoted by $V_{\rm
leaf}$.

If a tree is planted then there is a unique orientation
towards the root and hence each vertex has incoming edges
and at most one outgoing edge, the root being the exception in
having only incoming edges.
We will sometimes also use the
{\em arity}  $\arity{v}$ of $v$ to denote the number of incoming edges
to the vertex $v$.
Notice that for all vertices {\em except the root} $\val{1}=\arity{v}+1$,
but for the root $\val{v_{\rm root}}=\arity{v_{\rm root}}$.
In the figures the orientation of the edges toward the root is
taken to be downward.

For a tree $\t$ and $e\in E(\t)$ we will denote
the tree $\t'$ obtained from $\t$ by contracting $e$ by $\t'=\t/e$.
If in a rooted tree the marked flag $f_0$ is contracted, we fix the new
marked flag to be the image of the flag $f_{1}=N(f_0)$.
In this situation we will also say that $\t$ is obtained from $\t'$
by inserting an edge, and if we want to be more specific we
might add ``into the vertex $v$'', where $v$ is the image of $e$ under
the contraction and write: $e\mapsto v$.

If there is a vertex $v$ of valence $2$ in a tree, we denote
by $\t/v$ the tree $\t/e$ where $e$ is either
one of the two edges incident to $v$.
This just removes $v$ and splices together its two edges.

A {\em black and white (b/w) tree} is a pair $(\t,\color)$, that is
is a planar planted tree $\t$
whose set of vertices comes equipped with a map called color
$\color:V(\t)\to \Z/2\Z$, which satisfies that {\em all leaves are mapped
to $1$ and the root is mapped to $0$}.

 We call
the inverse images of $0$ black vertices and the inverse images of $1$ white
vertices.  The sets of black and white vertices will be
denoted by $\Vblack$ and $\Vwhite$ respectively. In particular, the
condition above then means that {\em  all leaves are white and the root
is black}.

 In a b/w tree the edges which have
two black vertices will be called {\em black edges} and denoted by
$\Eblack$. Similarly $\Ewhite$ denotes the white edges, that is
those whose vertices are both white. All other edges will be called
{\em mixed} and denoted by
 $\Emix$.
When contracting an edge, we fix that the color of the new vertex is black
if the edge was black and white if the edge was white. In the case
that the edge is mixed, we fix the color of the new vertex to be white.

A b/w tree is called {\em bipartite} if all edges are mixed. A b/w
tree is called {\em stable} if there is no black vertex $v_b$ with
arity $1$, except for the root which is the only black
vertex that may have valence $1$ and it may only have
valence $1$ if its unique incident edge is mixed.

A b/w tree is called {\em stably bipartite} if the following
conditions hold
\begin{enumerate}
\item There are no white edges.
\item There are no black vertices of arity $1$ and valence
$2$
both of whose incident edges  are black.
\item There are no black vertices of arity $1$ and valence $2$
where one edge is black and
the other edge is a leaf edge.
\item the root may have valence $1$, but only if its unique incident edge is mixed.
\end{enumerate}
Notice that a stably bipartite tree becomes bipartite when all
the black edges are contracted and stable if all the black vertices
of valence $2$ are removed. Stable trees and stably bipartite
trees are closed under contraction of black edges.

The effective white angles of a b/w tree are those angles whose
vertices are white {\em and}  have two distinct flags. They will be
denoted by $\wangles$. All effective white angles of flags at a given white
vertex $v$ will be denoted by $\wangles(v)$.

The conditions above are perhaps not so obvious from the tree point of
view but they are quite natural from an arc/foliation point of
view (see Appendix \ref{arcapp}).

We fix that a b/w subtree of a b/w tree has a white root.

An $S$--labeled b/w tree is a b/w tree  together with a bijective
labelling $\lab: S\to V_{\rm white}$; we will write
$v_i:=\lab^{-1}(i)$. When contracting a white edge, we label the new
white vertex by the union of the two labels considered as sets.

We will also need to cut and assemble a tree by gluing subtrees
along a tree. The basic operation is {\em replacing a vertex with a tree}. Combinatorially this
is defined as follows. {\em Replacing a black vertex} $v$ in a planar b/w  tree $\t$
by a planar b/w tree $\t'$
whose number of leaves equals $\arity{v}$ and whose root is black
 means the following: (1)
 we remove all flags incoming to $v$ from $\t$;
(2) we add all the vertices of $\t'$
that are not leaves and all flags of $\t'$
except the flags incident to the leaves;
(3) since the cardinality
of the sets of flags incident to leaves of $\t'$ and
the set of incoming flags are the same and both of them have an
order, there is a unique bijection $\phi$ preserving this order.
We  ``glue in'' the new vertices and flags by keeping $\imath$ wherever
it is still defined and using $\phi$ and $\phi^{-1}$ for
the other flags. We also fix that the outgoing flag of $v$
has the root of $\t'$ as its vertex.

When {\em replacing a white vertex} $v$ of a planar b/w tree $\t$ by a
planar b/w tree
we proceed as follows: (1) we remove the vertex $v$
and all incoming flags of $v$ from $\t$; (2) we add all the vertices
of $\t'$ that are not leaves and all flags of $\t'$ which are not
incident to the leaves; and 
(3) we glue the flags as in the case of replacing a black vertex.
There is a special case, in which 
a white vertex that is adjacent to a root of valence one is replaced.
In this case, as a final step, we contract the unique edge
incident to the root.

See Figure \ref{treecut} for an example. The example has
extra labellings, which is discussed in \S\ref{treecutsection}.

We will deal with three sets of trees in particular:

\begin{df}
We define $\bipartite(n)$ to be the set of  $\{1,\dots,n\}$--labeled
b/w bipartite planar planted trees. We use $\bipartite$ for the collection
$\{\bipartite(n),n\in {\mathbb N}\}$.
\end{df}
We recall that we fixed that all leaves of a b/w tree are
white and the root is black.

\begin{df}
We let $\stable(n)$ be the set of  $\{1,\dots,n\}$--labeled b/w
stable planar planted trees. We denote by $\stable$ the collection $\{\stable(n),n\in {\mathbb N}\}$.
\end{df}

\begin{df}
We let $\heights(n)$ be the set of pairs $(\t,\height)$, where
$\t$ is a black and white stably bipartite trees  $\{1,\dots,n\}$--labelled
and
$\height:E_{\rm black} (\t)\to \{\var, 1\}$ called the height function.
The collection $\{\heights(n),n\in {\mathbb N}\}$ will be simply be denoted by $\heights$.
\end{df}
 Here $\var $  stands for variable height.
We will denote the set of edges labelled by $\var$  by $\varedges$.

\begin{nota}
 We will use the notation $\Z S$, for the Abelian group generated by  a set $S$.
E.g.\ $\Z\bipartite(n)$ and $\Z\stable=\bigoplus \Z\stable(n)=\M$.
\end{nota}

\subsection{The differentials}

There are natural differentials on each of the three Abelian groups $\Z\bipartite$, $\Z\stable$ and
$\Z\heights$. The differential for $\bipartite$ was given in \cite{del} and the one on
$\Z\stable$ was introduced in \cite{KS}. We will briefly recall the definitions and give a new definition
for a differential on $\Z\heights$.

\subsubsection{The differential for $\Z\bipartite$}
Following \cite{del, cact},
we fix a tree $\t\in \bipartite(n)$ for each effective
white angle $\a\in \wangles$ and let $\del(\a)(\t)$ be the
tree obtained by {\em collapsing the angle $\a$}.
Combinatorially put let $\a=\{f_1,f_2\}$, let $e_i=\{f_i,\imath(f_i)\}$
and set $v_i=\del\imath(f_i)$. Then  $\del(\a)(\t)$ is the tree
where $v_1$ and $v_2$ are identified as are $e_1$ and $e_2$.
The new marked flag will simply be the image of the original marked flag
(see Figure \ref{collapse}). Using this notation,
the differential is defined as:

\begin{figure}
\epsfxsize = .8\textwidth \epsfbox{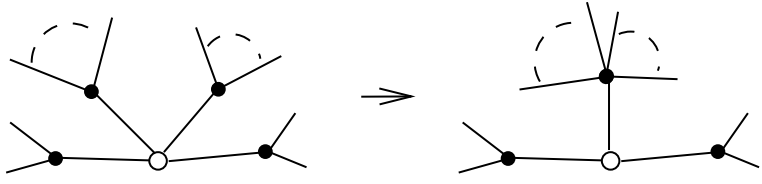}
\caption{\label{collapse} Collapsing a white angle.}
\end{figure}

\begin{equation}
\del(\t)=\sum_{\a\in \wangles}\pm\del(\a)(\t)
\end{equation}

\subsubsection{The differential on $\M=\Z\stable$}
Following \cite{KS},
we fix a tree $\t \in \stable(n)$. We will consider all trees that are obtained from $\t$ by adding an
edge which is either mixed or black.
 That is, the summands of the differential are indexed by pairs $(\t',e)$
such that the tree $\t'/e$ obtained by
contracting $e$ is equal to $\t$ and $e\in \Eblack\amalg \Emix$.
Here the cyclic structure is the induced one and we recall that
 the rules for contracting edges prescribe that
the image of a black edge is a black vertex
and the image of a mixed edge is a white vertex.

\begin{equation}
\del(\t)=\sum_{\small\begin{array}{c}(\t',e)\\
\t'/e=\t,e\in \Eblack\amalg \Emix\end{array}}  \pm \t'
\end{equation}

Alternatively one can sum over local contributions $\del(v)(\tau)$ considering
only those edges whose image is $v$. This is the way it is written in \cite{KS}.
\begin{equation}
\del(\t)=\sum_{\small\begin{array}{c} v\in \Vwhite, (\t',e)\\
\t'/e=\t,e\in \Emix, e\mapsto v\end{array}}\hspace{-.5cm} \pm \t'
+\sum_{\small\begin{array}{c} v\in \Vblack, (\t',e)\\
\t'/e=\t,e\in \Eblack, e\mapsto v\end{array}}  \hspace{-.5cm} \pm \t'
\end{equation}

\subsubsection{The differential on $\Z\heights$}

We now fix $(\t,\height)\in \heights(n)$. For the differential, we will sum
\begin{itemize}
\item[a)] over collapsing the white angles, i.e.\ elements of $\wangles$
and
\item[b)] over contracting the black edges labeled by
$\var$.
\end{itemize}

For a white angle $\a\in \wangles$,
we again let $\del_{\a}(\t)$ be the tree with the white angle collapsed.
We can keep the height function since the collapsing angles
does not affect the set of black edges --- only two mixed edges are identified.
 For an edge $e\in \varedges\subset E_{\rm black}$ we set
 $\del_{e}(\t,\height)=
(\t/e,\height|_{\Eblack \setminus e})-(\t,\height')$ where $\height'(e)=1$
 and $\height'(e')=\height(e')$ for $e\neq e'$. The differential is now
 \begin{equation}
 \del(\t)=\sum_{\a\in \wangles}\pm\del_{\a}(\t)+
\sum_{e\in \varedges} \pm \del_{e}(\t)
  \end{equation}

\subsubsection{Signs}
\label{signsection}
As usual the fixing of sign conventions is bothersome,
 but necessary. The quickest way is to use
 tensor products of lines of various degrees indexed by
the sets of edges and/or angles. See \cite{cact,del,KS} for detailed
discussions. One way to fix an order of the tensor factors is to fix an
enumeration of all flags by going around the planar planted tree
starting at the marked flag and then using the map $\imath$ and the
cyclic order to enumerate.  Hence all vertices, the subset of white vertices,
angles, the subset of white angles, and edges are enumerated by counting them when their first
flag appears. We will call this the {\em planar order}.
To fix the signs one simply fixes weights of the elements of the
ordered sets.

A third way, and perhaps the cleanest for the present discussion, is to use the
geometric boundary of polytopes as we will discuss in \S \ref{cellsection}
below. In particular, the signs for the different types
of trees are fixed by equations (\ref{bipartcelleq}),
(\ref{stablecelleq}) and (\ref{heightcelleq}).

  \begin{prop}
  In all three cases $\bipartite,\stable,\heights$ the map $\del$ satisfies $\del^2=0$.
  \end{prop}

  \begin{proof}
  In all cases this is a straightforward calculation.
The signs are such that inserting two edges or alternatively collapsing two
  edges or angles (or one edge and one angle) in different
orders yields the same tree, but with opposite signs, since
these elements are ordered
and formally of odd
degree in any of the above formalisms.
  \end{proof}

\subsubsection{The maps $\stinf$ and $\ininf$}
There are maps $\stinf:\Z\stable \to \Z\bipartite$ and
$\ininf:\Z\bipartite \to \Z\stable$ which were defined in \cite{del}.

The first map $\stinf$  is given as follows. If there is a black vertex
of valence $> 3$, then the image is set to be $0$. If all black
vertices are of valence $3$, we (1)  contract all black edges and (2)
insert a black vertex into each white edge, to make the tree bipartite.
It is clear that the leaves will stay white.
The global marking, viz.\ root is defined to be the image of the marking
under the contraction.

The second map $\ininf$ is given as follows: (1) Remove all black
vertices whose valence equals 2 and (2) replace each black vertex of valence
$>2$ by the binary tree, with all branches to the left.
This is of course not symmetric, but any choice will do.
Now we see that $\stinf$ is surjective, since
$\stinf\circ\ininf=id$.

\begin{lem}\cite{del}
\label{stinflem}
The  map $\stinf$  behaves well with respect to the differential.
$\stinf(\del(\t))=\del \stinf(\t)$.
\end{lem}

\begin{proof}
This is a straightforward calculation, see \cite{del}.
\end{proof}

\subsection{Operad structures on $\Z\bipartite$ and $\M=\Z\stable$}

Both the operad structures are what one could call an insertion
operad structure. They have been previously defined in \cite{del} and in
\cite{KS} respectively. The latter was defined
combinatorially in \cite{KS}, but also can be induced from the topological
level; see Appendix \ref{arcapp} and Proposition \ref{textpropone}.

There are two equivalent ways to describe this type of operation.
The indexing is always over the white vertices. 
Inserting a tree $\t'$ into a tree $\t$ at the vertex $v_i$ means the signed sum over all trees $\t''$
which contain $\t'$ as a sub--tree such that $\t''/\t'=\t$ with the
image of $\t'$ being $v_i$.

\begin{equation}
\t\circ_i \t'=\sum_{\t'':\t''/\t'=\t, \t'\mapsto v_i}\pm \t''
\end{equation}

Here one also fixes that $\t''$ be either in $\bipartite$ or $\stable$.
Also, contracting $\t'$ as a subtree in the case of $\bipartite$ means that
we first insert an additional black  edge for the black root of the subtree,
such that the new vertex has valence $1$ when considered as a  vertex of the
subtree,
 and then contract the subtree; the result would not
be bipartite otherwise.
In the case of a stable tree, there is the provision that if the root
of $\t'$ has valence $1$ then the root edge is contracted before
identifying $\t'$ as a subtree, i.e.\ this vertex is not present in the
subtree.
 The sign again is given by one of the schemes in \S \ref{signsection}.

Alternatively, one can describe a 3--step procedure consisting of first cutting
off all the branches over $v_i$, then grafting $\t'$ into $v_i$, and finally
grafting the branches back to $\t'$ keeping their order as induced
by the cyclic order on $\t$.
We refer to \cite{KS,woods,cact,del} for more details.

\begin{prop}\cite{KS,del}
The collections $\stable$ and $\bipartite$ are dg--operads.
\end{prop}
\qed

\begin{prop}\cite[Proposition 1.5.8]{del}
$\stinf$ is a morphism of dg--operads.
\end{prop}
\qed

\subsubsection{Operad structure on $\Z\heights$}
Strictly speaking, we will not need an operad structure on
$\Z\heights$ to prove the Main Theorem. However, there is indeed an
operad structure, and it and the operad structure on $\Z\bipartite$
can be understood as special cases of a operad structure  induced by
the quasi--PROP structure of Sullivan chord diagrams of \cite{hoch1};
see Appendix \ref{arcapp}.

\begin{figure}
\epsfxsize = .8\textwidth \epsfbox{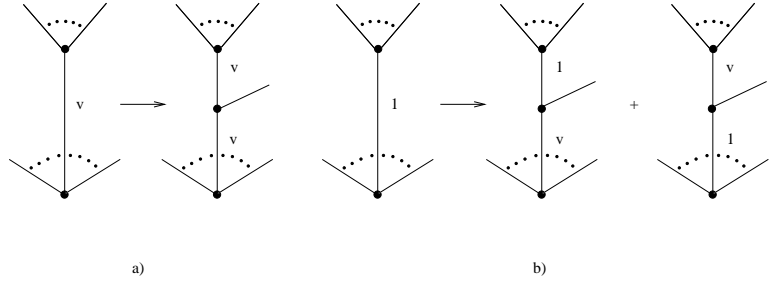}
\caption{\label{glueintoedge} Gluing a branch into an edge a) in $\Evar$ and b) of height $1$}
\end{figure}

We first give the definition combinatorially.
Given $(\t,\height)$ and $(\t',\height')$ we define $\mathcal S$
to be the following set of trees with height
$(\t'',\height'')$. $\t''$ is obtained by cutting the branches
of $\t$ above $i$, gluing in $\t''$ at $i$ and then
gluing in the branches in their planar order to the white angles
of the image of $\t''$ and into the black edges $\Eblack(\t'')$.
To glue a branch into an edge, we add a vertex to the edge
and glue the branch to this new vertex.
The admissible height functions $\height''$ coincide
with the original
height functions on all images of edges of $\Eblack(\t')$ and all
unaffected edges of $\Eblack(\t'')$. Let $e$ be a black
edge that has been split into $n$ black edges by gluing in $n-1$ branches.
If $e\in \Evar$ then {\em all} the values of $\height''$
on the edges it is split into are $\var$. If $e$ is labelled by $1$
then all but one of the labels are $\var$ and one label is $1$.
All these labels are allowed; see Figure \ref{glueintoedge}

\begin{equation}
(\t,\height)\circ_i(\t',\height')=\sum_{(\t'',\height'')\in {\mathcal S}}
\pm (\t'',\height'')
\end{equation}

\begin{prop}
The collection $\Z\heights$ yields a dg--operad.
\end{prop}
\begin{proof}
Somewhat tedious but straightforward calculation; or see Proposition \ref{appheightprop}
\end{proof}

\subsubsection{$\Ainf$ algebras}

Notice that the trees $\plantree$ in $\stable$ with $V_{\rm white}=V_{\rm leaf}$ form a
sub--operad $\Z\plantree$ of $\Z\stable$.

It is straightforward to see
that this operad is isomorphic to the operad of planar planted
trees with labelled leaves with the operation of grafting at the leaves.
Keeping this in mind
the following definition goes back to Stasheff (see \cite{MSS} for a
more complete history):
\begin{df}
An $\Ainf$ algebra is an algebra over the dg--sub--operad $\plantree$.
\end{df}

In particular, on an $\Ainf$ algebra $A$
there is an $n$--ary operation $\mu_n$ for
each $n\in {\mathbb N}$, such that $\mu_1$ is a differential $\del$,
and $\mu_2$ is associative up to the homotopy $\del(\mu(3))$.
After this there is a whole tower of homotopies governed by the
combinatorial structure of the $K_n$.

\subsubsection{Associative algebras}

We can also consider $\Z\corollas$, that is the bipartite
trees with white leaves only, as a sub--operad of $\Z\bipartite$.

\begin{lem}
 $\Z\corollas$ is isomorphic to the operad for associative algebras.
\end{lem}

\section{Polytopes and Trees}
\label{polysection}
In this section, we review associahedra and cyclohedra emphasizing
that they together with the standard simplex can be thought of
as compactifications of the open simplex. This in turn has an interpretation
as a configuration space.

\subsection{Simplices}
\label{simplexsection}
We let $\Delta^n$ be the standard simplex and $\dot \Delta^n$ be its
interior.

\subsubsection{Configuration space interpretation}

If we realize the simplex as
 $\Delta^n = \{t_1,\dots, t_{n} | 0\leq t_1\leq \dots \leq
t_n\leq 1\}$ and  $\dot\Delta^n=\{t_1,\dots, t_{n} | 0< t_1<\dots< t_n<1\}$,
then $\dot \Delta^n$ is the configuration space of $n+1$ distinct
points on $\dot I=(0,1)$ and the closure just lets the points
collide with each other or with $0$ and $1$. That is, the space
is just the compactification obtained from $n$ unlabelled,
not necessarily distinct,
points on $[0,1]$.

The interior of this compactification
 is the same as considering $n$ distinct points on
$S^1$ with one point fixed at $0$.  The compactification then
distinguishes if the points collide from the right or left with
$0$, but keeps no other information.

\subsubsection{Tree interpretation}
As a polytope, the simplex is a CW complex and of course the cells are 
again just simplices. We can give a tree interpretation as follows:
the cell defined by an $n$ simplex will be indexed by a tree
$\starwhite_n$ which we
call a white star. The tree $\starwhite_n$ is the unique bi--partite
tree with exactly one white
vertex that is not a leaf, of which there are $n$, and all of whose non--root
black vertices have valence $2$ and the root has valence $1$.
We can pictorially think of the white vertex as $S^1$ and the incident
 edges as indicating the points on $S^1$, where the root marks $0$.
The boundary map is just
the sum of collapsing of the white angles.
After collapsing an angle, we still have only one white non--leaf
vertex, but the black vertices may have valence $2$ or the root may have
valence $3$. The leaves incident to a black non--root vertex are
the points that have collided with each other and the leaves
incident to the root are the points that collided with $0$. Since
the tree is planar, we can distinguish if this happened from the right or
left.

\subsubsection{Topological interpretation}
\label{topheightsection}
We can make the cell decomposition above topological as follows.
To each white angle of $\starwhite_n$ we associate a number in $(0,1]$ that
is we have a map $\topheight:\wangles(\starwhite_n)\to (0,1]$,
which we subject to the
condition that the total angle
at the white vertex is 1: $\sum_{\a\in \wangles}\topheight(\a)=1$.
If the only white angle is not effective, we can just label it by $1$.
We can imagine that these angles measure the distance between the points of $S^1$ in units of $2\pi$. The open part is then just $\dot \Delta^n$ and the
closure is $\Delta^n$. The boundary comes from sending the length of
the angles to zero and collapsing the angles.

\subsection{Associahedra}
The associahedra are abstract polytopes introduced by Stasheff
\cite{stasheffass1, stasheffass2} and fittingly are also called Stasheff polytopes. The vertices of the associahedron $K_n$  correspond to
the possible full bracketings of the expression $(a_1 \cdots a_n)$,
e.g. $(((a_1a_2)a_3)(a_4 a_5))$. Each such bracketing can be
depicted as a planar planted tree by thinking of the bracketing as giving
a flow chart. The dimension $l$ faces correspond to bracketings
which are missing $l$ pairs of brackets; here it
is assumed that the outside bracketing is always present. The highest dimension
and hence the dimension of $K_n$ is $n-2$. We will also allow and use $K_2=pt$.
 E.g.\ $((a_1a_2)a_3a_4)$ is of dimension $1$ and $(a_1a_2a_3a_4)$ is
of dimension
 $2$. The boundary of the faces is given by inserting one set of brackets in all possible ways.
In the tree picture the codimension is given by the number of internal,
that is non--leaf edges
and the boundary map is defined by inserting an edge in all possible ways.
It is a well known fact that the faces of $K_n$ are products
$K_i\times K_{n-i}$.

\subsubsection{Labelling} It will be
convenient to use other indexing sets and consider $S$--labelled
associahedra $K_S$. In the bracket formalism this is the indexing set of
the elements $a_i$. This is already useful in the description
of the boundary, since the boundary components are distinguished
by different labels. In particular the boundary is given
by
\begin{equation}
\del K_n=\sum_{(I',II'')}K_{I'}\times K_{I''}
\end{equation}
where $I'=\{j,\dots j+k\}$ with $1\leq j, k\geq 1, j+k\leq n$ and
$I''=\{1,\dots,j-1, I', j+k+1,\dots n\}$. This choice corresponds
to the bracketings compatible with $(a_1\cdots a_{j-1}(a_j\cdots a_{j+k})a_{j+k+1}\cdots a_n)$.

\subsubsection{Configuration interpretation}
The space $K_n$ can be viewed as a ``real Fulton--MacPherson
compactification \cite{FM}'' of the space of $n-2$ distinct points on the
interval (0,1) \cite{MSS}. The information that is kept are the relative
speeds of multiple collisions. Just as above by identifying
$0$ and $1$ one can view
this as a compactification of distinct points on $S^1$, where now
one keeps track separately of the points colliding with $0$ from the right and from the
left and of the relative speeds of these two processes.

\subsubsection{A first CW realization with stable trees $\Tass$}
As an abstract polytope the associahedra are naturally CW complexes. The
cells for $K_n$ are indexed by planar planted trees with $n$ leaves
and their dimension is given by $n-2-|E|$. We will make
the leaves white and consider them to live in
$\stable(n)$ and insist that the labelling from $1$ to $n$ is
consistent with their planar order.
To be precise we let $\Tass(n)$ be the trees in $\stable(n)$
whose only  white
vertices are leaves.
Each cell $\cell(\t)$ represented
by a tree $\t\in \Tass(n)$ is a product

\begin{equation}
\cell(\t)=\times_{v\in V(\t)} K_{\arity{v}}
\end{equation}

\begin{figure}
\epsfxsize = .8\textwidth \epsfbox{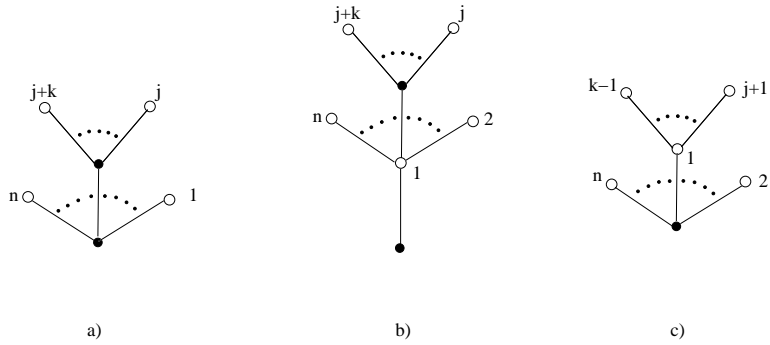}
\caption{\label{boundaryW} a) Boundary trees of $K_n$, b),c)
Boundary trees of $W_n$}
\end{figure}

The differential given by taking the boundary agrees with the sum over
all possibilities of inserting a black edge
which is  the one inherited from $\stable$, i.e.\
$\del(\cell(\t))=\cell(\del(\t))$, where we extend $\cell$ in
the obvious fashion to linear combinations. Notice that the labelling
sets are now induced by contracting either all  the edges of the ``upper''
vertex or alternatively contracting
all the edges of the ``lower'' vertex; see Figure \ref{boundaryW}a).

\subsubsection{A second CW realization with trees with heights $\Tassheight$}
\label{secondasspar}
There is an alternate natural CW structure which is actually a cubical
decomposition of the associahedra.
This is sometimes called the Boardman--Vogt decomposition \cite{BV} where
strictly speaking it is a Boardman--Vogt construction for the
operad of monoids; see also \cite{MSS,KS}.
The cells of this compactification are cubes and are indexed by particular
trees in $\heights$. The trees are those in which all the white
vertices are leaves, viz.\ $\Tass(n)$ and again we insist that
their planar order is given by the labelling.
Putting  all possible height functions
on them, we obtain a subset $\Tassheight(n)\subset \heights(n)$.
The cell indexed by $\t$
is
\begin{equation}
\cell(\t)=I^{\Evar}=\times_{e\in \Evar}I
\end{equation}

The boundary is  given
by using the differential of $\heights$. We again have that
$\del(\cell(\t))=\cell(\del(\t))$, where we extend $\cell$ in
the obvious fashion to linear combinations.

\begin{rmk}
Notice that this CW decomposition is a subdivision of the first.
The cells of the finer decompositions that belong to a given cell given by
a tree $\t$ can be described as follows: first label all black edges of $\t$ by $1$ and then consider all trees in $\heights$ which can be contracted
to $\t$ and whose labels match on the non--contracted edges.
\end{rmk}

\begin{rmk}
We actually rediscovered this decomposition from the arc point of
view; see Appendix \ref{arcapp}. After presenting the results, we
realized that this decomposition coincides with
a Boardman--Vogt construction, but we would like to point
out that it also comes naturally from
a topological
quasi--operad; see Appendix \ref{arcapp}
\end{rmk}

\subsubsection{A topological realization via trees with heights}
Since their introduction, people have looked for convex
polytope realizations of the associahedra.
This has lead to several nice results and constructions; see
e.g.\  \cite{CD,CFZ,FR,L} for recent results and
also \cite{MSS} for more references and details.

Taking the cue from the above cell decomposition
one can easily give a realization which is {\em not a convex
polytope}, but a {\em PL} realization.
For this we will consider the trees with bounded heights,
that is pairs $(\t,\topheight)$ where $\t\in \Tassheight$
and $\topheight:E_{\rm black}\to (0,1]$. If we let $\E(n)$
be the set of all possible black edges for such trees with fixed $n$,
this space is naturally a {\em subspace} on
$\mathbb{R}^{\E(n)}$.

Notice that in the subspace topology the limit where $\height(e)\to
0$ for some edge $e$ is naturally identified with the tree with heights, where
this edge has been contracted. Moreover the boundaries are  also
naturally given by the same PL realization.

\begin{prop}
The construction above yields  a PL realization of $K_n$.
\end{prop}
\qed

\begin{figure}
\epsfxsize = \textwidth \epsfbox{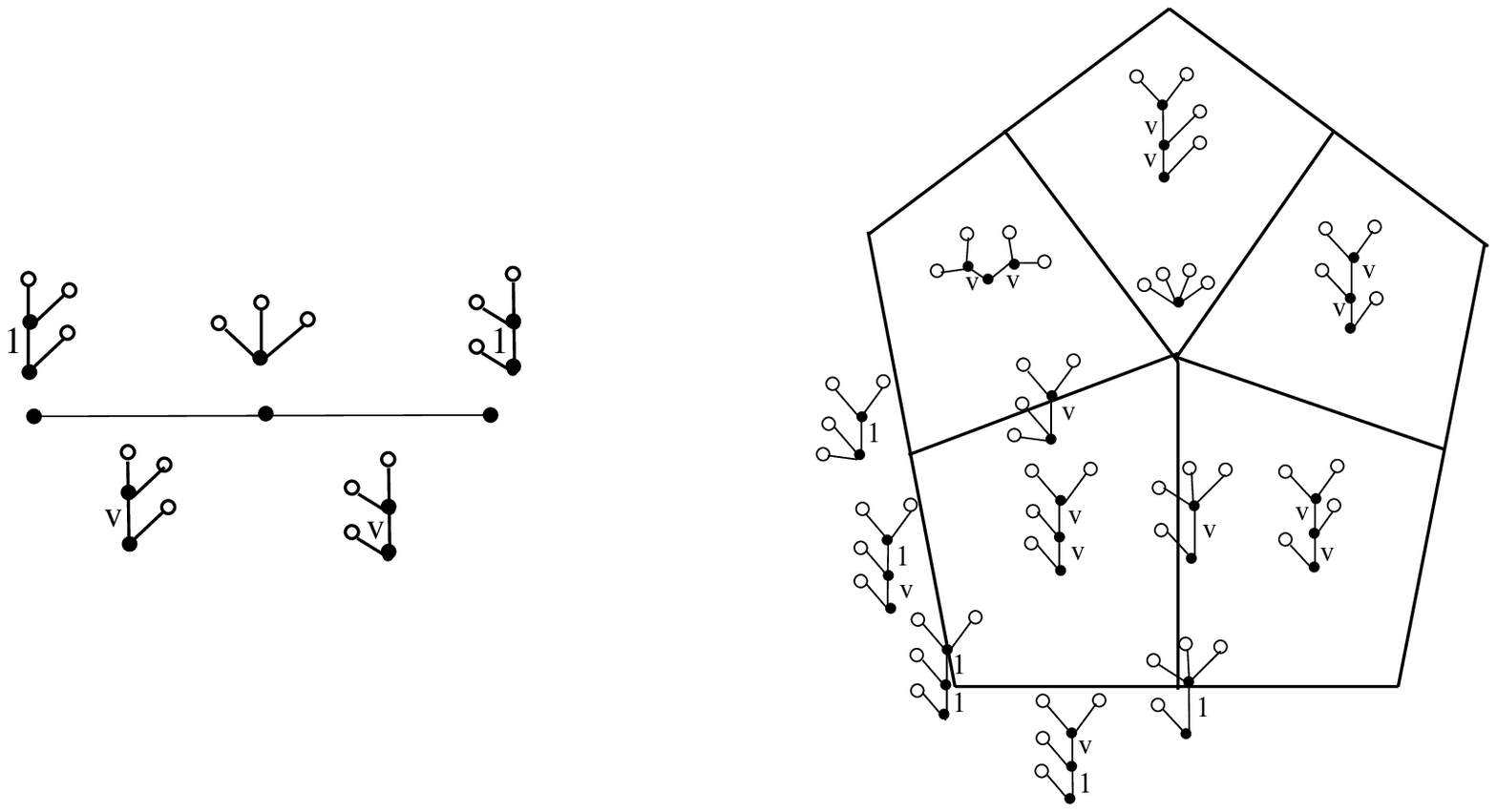}
\caption{\label{k4decomp} The decompositions of $K_3$ and $K_4$.
For $K_4$ the trees of dimension less
than two are only given for the lower--left cell}
\end{figure}

\begin{df}
We call a topological height function $\topheight$ on a tree with heights
$(\t,\height)$ compatible if $\topheight(e)=1$ when $\height(e)=1$
and $\topheight(e)\in (0,1)$ when $\height(e)=\var$.
\end{df}

The elements inside a given cell $(\t,\height)$ are then the elements
 $(\t,\topheight)$ with $\t\in \Tass$ and $\topheight$ a compatible
topological height function. The elements in the closure of this set,
that is those on the boundary of the cell
are those pairs ($\t',\topheight)$  where $\t'$ can
be obtained from $\t$ by contracting any number of edges of $\Evar$,
 $\topheight$ may now take values in $(0,1]$,
and at least one edge is contracted or one edge
$e\in \Evar$ has $\topheight(e)=1$.

\subsection{Cyclohedra}

The cyclohedra are abstract polytopes introduced by Bott and Taubes
\cite{BT}. The vertices of the cyclohedron $W_n$ correspond to the
possible full cyclic bracketings of the expression $a_1 \cdots a_n$,
e.g. $a_1))(a_2(a_3$. The $l$ dimensional sides are given by the
bracketings missing $l$ brackets. Here we allow the empty
bracketing. The boundary map is given by inserting one pair of
brackets in all possible ways. The dimension of $W_n$ is $n-1$. We
will also allow and use $W_1=pt$. Moreover, as with the $K_n$,
we will need to consider $S$--labelled $W_n$, that is $W_S$, where
$S$ is the indexing set of the elements.

It is well known and easy to check in this formalism, that the
$\codim(l)$ cells are products of $l$ polytopes of which one is a
cyclohedron and the others are associahedra.
 The possible sub--bracketings of a cyclic bracketing are
given by independent choices of regular bracketings.

From the description above, we see that the boundary is given by
\begin{equation}
\label{boundarycyceq}
\del(W_n)=\sum_{(I',I'')} W_{I'}\times K_{I''}
\end{equation}

Here the indexing sets on the right hand side are the ordered sets
$I''=(j,j+1,\dots,j+k)$ $j\leq 1, j+k\leq n$ for $k\geq 1$ and
$I'=(1,\dots, j-1,I'',j+k+1,\dots,n)$, or $I''=(2,\dots
j,\{1,j+1,\dots,k-1\},k,k+1,\dots, n)$ for $j<k$  and $I'=(\{1\}\cup
I'',j+1,\dots, k-1)$, here if $k+1=j$, means that $I'=(1)$.

Again these indexing sets follow from contracting the
relevant edges of the
``upper'' or ``lower'' vertex, see Figure \ref{boundaryW} b), c).

\subsubsection{A configuration interpretation}
The way they were originally introduced by Bott and Taubes
they are the blow--up of a configuration space.
This is also related to the Axelrod--Singer \cite{AS} compactification
of configuration space, see \cite{MSS} for details.
In particular the cyclohedron  $W_n$
is the compactification of the configuration of $n$ distinct
points on $S^1$ with one point fixed at $0$, see \cite{MSS} for details.

\subsubsection{A first CW realization in terms of stable trees $\Tcyc$}
Again, we have the natural structure of CW complex. A tree depiction
is given as follows: we consider trees which are
 $\{1,\dots,n\}$--labelled b/w stably bipartite
with at most one white internal vertex labelled by $1$ and
 all other white vertices are leaves and these leaves
 are labelled commensurate with the planar order. This means
 that if there is an internal white vertex, all the leaves
are labelled $2,\dots,n$ in that order and if there is no internal
white vertex all white vertices are leaves and the order of the
leaves labelled $2,\dots, n$ is exactly this order, while the vertex
labelled by $1$ may appear anywhere in the planar order. We will
call these trees $\Tcyc$.
 The ``big''
cell representing the whole cyclohedron is the unique tree which has
no black vertices. Again, we can think of the internal white vertex
as $S^1$ and its edges as indicating the location of the points, if
we wish.

The boundary comes from inserting a mixed edge into the white non--leaf vertex,
which yields a product of a cyclohedron and an associahedron.

In general we have that
the cell of $\t$ is given by
\begin{equation}
\cell(\t)=\times_{v\in \Vwhite} W_{\val{v}}
\times \times_{v\in \Vblack} K_{\arity{v}}
\end{equation}

 The differential
is then the differential of $\stable$,
$\del(\cell(\t))=\cell(\del(\t))$, where we extend $\cell$ in the
obvious fashion to linear combinations.

\subsubsection{A second CW realization in terms of trees with heights}
\label{secondcycpar}

We will exhibit another CW realization for $W_n$ which has the
following trees as an indexing set: these are the trees
$\Tcycheight$ in $\heights$ which have  $n$ white vertices and at
most one white non--leaf vertex. We consider these trees to be
labelled by $\{1,\dots,n\}$ and impose the same conditions as for
$\Tcyc$, i.e.\ the vertices $v_2,\dots,v_n$ are leaves and the
planar order of this subset is the one written. The vertex $v_1$ may
be internal and may appear anywhere in the planar order of all white
vertices, even if it is a leaf.

 We define a cell of such a tree as
\begin{equation}
\cell(\t)=\times_{v\in \Vwhite} \Delta^{\arity{v}} \times I^{\Evar}
\end{equation}

We now get a CW complex $\CWcyc(n)$ by using the trees above and
 the differential of $\heights$ to
define the boundary and hence the attaching maps.

To fix terminology we will call a black vertex {\em potentially unstable}
if it is adjacent to a non--leaf mixed edge.

\begin{lem}
\label{factslem}
The following statements hold for the CW complex $\CWcyc(n)$

\begin{itemize}

\item [(i)] The dimension of $|\CWcyc(n)|$ is $n-1$.
 The top--dimensional cells are precisely indexed by
the trees such that there are only $n-1$ leaves,
the arity  of all black vertices
is $\leq 2$,
 all potentially unstable
non--root vertices are valence $2$, the root is either
 not potentially unstable or if it is, it is of arity $1$,
{\em and}
all black edges are labelled by $\var$.

\item[(ii)] All $0$--cells are indexed by trees whose white
vertices are all leaves, and all black edges  have height $\height$
equal to one.

\item[(iii)] All $k$--cells are
in the boundary of $k+1$ cells for $k<n-1$ and each chain of cells
such that the successor is in the boundary of the predecessor  has
length $n$.
\item[(iv)] The codimension $1$ cells are given by trees of the following
types:

\begin{itemize}
\item[(a)] A tree as in (i) with only one black edge labelled by $1$
\item[(b)] A tree as in (i) but with one of the non--root potentially unstable
vertices having valence $3$.
\item[(c)] A tree as in (i) but with one of the other black vertices
(not potentially unstable) of valence $4$.
\item[(d)] A tree as in (i) but the root vertex not potentially unstable having valence $3$.
\item[(e)] A tree as in (i) but the root vertex potentially unstable and of valence $2$.
\item [(f)]  A tree as in (i) but no internal white vertex.
\end{itemize}

Each cell of the types (b), (c), (d)  and (e) are in the boundary of precisely
2 top--dimensional cells and the cells of type (a) and (f) are in the boundary
of exactly one top--dimensional cell.
\end{itemize}
\end{lem}
\begin{proof}
Ad (i), by counting dimensions, we see that the dimension of
cells listed is indeed $n-1$. It is also just a dimension count
that these cells are indeed the maximal ones. Any higher arity of a
black vertex or a black
edge labelled by anything else but $\var$ will lead to a dimension drop as one
could change the label to $\var$, insert a new edge, or ``split'' an angle.

This procedure also shows the claim (ii) and (iii). The chains are given
by a series of a total number of $n-1$ contractions and collapsing.

To be in codimension $1$ the dimension count has to go down by one from
the top--dimensional cells by moving to the boundary.
Starting with a top--dimensional cell indexed by a tree with heights,
we can (1) relabel an edge from $\var$ to $1$,  (2) contract an
edge labeled by $\var$ or (3) collapse one white angle. The result of (1)
will be a tree of type (a), the result of (2) will be of type (b) if
the edge was incident to a potentially unstable vertex and of type (c) if
it was not and not incident to the root. It will be of type
(d) if it was adjacent to the root and after contraction
the root is not potentially unstable. It is of type (e) if the root
becomes potentially unstable.

The results of (3) will be of type (b) if the angle did not
have the root as one of its vertices and will be of type  (e) or (f)  if it did.
This may only occur if the root had valence $1$.

To determine the cells that lead to the particular boundary, we reverse
the above operations in all possible ways.
In  case (a) we can only re-label
the edge by $\var$ and in case (f) the only possibility is
to ``split'' the angle of the vertex labelled by $1$
at the root in order to obtain a non--leaf
white vertex.

In case (b) the only two possibilities are to insert
a black edge labelled $\var$
or to ``split'' the vertex into a white angle.
In  case (c) there are exactly two different ways to insert
one black edge labelled by $\var$, this is analogous to
the case of $K_3$. The case (d) is analogous. Finally,
in the case (e) we can either insert an edge marked $\var$
to make the root not potentially unstable, or split the angle.
\end{proof}

\begin{thm}
\label{cycthm}
The CW complex $\CWcyc(n)$ is a CW  realization of the cyclohedron.
This is a refinement of the polytope CW complex. The additional
$0$--cells correspond to the refinement of the associahedra.
\end{thm}
\begin{proof}
We will make the proof by induction. We have to show that the
boundary of $\CWcyc(n)$ is indeed composed of $W_{n-i}\times
K_{i}$'s with $i\geq 2$. First the case of $n=1,2$ are trivial to
check. Here we use a decomposition of these polytopes viewed as cell
complexes known by induction for the cyclohedra and the previous
results for the associahedra. The case $n=3$ is in Figure \ref{W3decomp},
and the case of $n=4$ is worked out in Appendix \ref{cycloapp}.
We let $\omega(n)=\sum_{\t: \dim(C(\t))=n-1}\t$ be the sum of all
top--dimensional cells. Now $\del\omega=\sum \del\t$ and on the
right hand side we will only have the terms of the types (a) and (f)
of the lemma above, since the terms of type (b)--(e) cancel out.
For terms of type (f) we notice that they sum up to associahedra
$K_{n}$, labelled by the different orders of $1,\dots, n$ which
respect the natural the order of $2,\dots, n$. I.e.\ all the faces
of the cyclohedron which are associahedra, using the second CW
decomposition described above. For terms of type (a) we first notice
that the cells are products of the cells associated to the trees above and below the black
edge marked by one. To be precise given a tree $\t$ of the type (a)
with the edge $e$ marked by $1$ we let $\t'$ be the tree with $e$
and  all the edges above $e$ contracted and  $\t''$ be the subtree
of $\t$ above $e$. Then the cell $\cell(\t)=\cell(\t')\times
\cell(\t'')$. The cell $\cell(\t'')$ has no internal white vertex
and is part of an associahedron. The cell $\cell(\t)$ has a white
vertex and by induction this is part of lower dimensional $W_k$.
Fixing either tree, i.e.~$\t'$ or $\t''$ and regarding all the possible
trees they can come from,
we see that the summands needed
to complete the associahedron , as discussed
in \S \ref{secondasspar}, and the cyclohedron, as in the assumption, which
we have established per
induction for the boundary terms of lower dimension, are all realized.
Moreover it is straightforward to check that all the needed labellings enumerated
in equation (\ref{boundarycyceq}) are realized and only those.
By Lemma \ref{factslem} the CW complex made up out of the
consolidated cells then yields an abstract polytope
 and this polytope is the
 cyclohedron $W_n$.

Finally, the $0$--cells are indexed by trees with no effective
white angles and hence all white vertices are leaves. All
the black edges are labelled by $1$ and hence these correspond
exactly to the $0$--cells of the respective associahedra.
\end{proof}

\begin{figure}
\epsfxsize = \textwidth \epsfbox{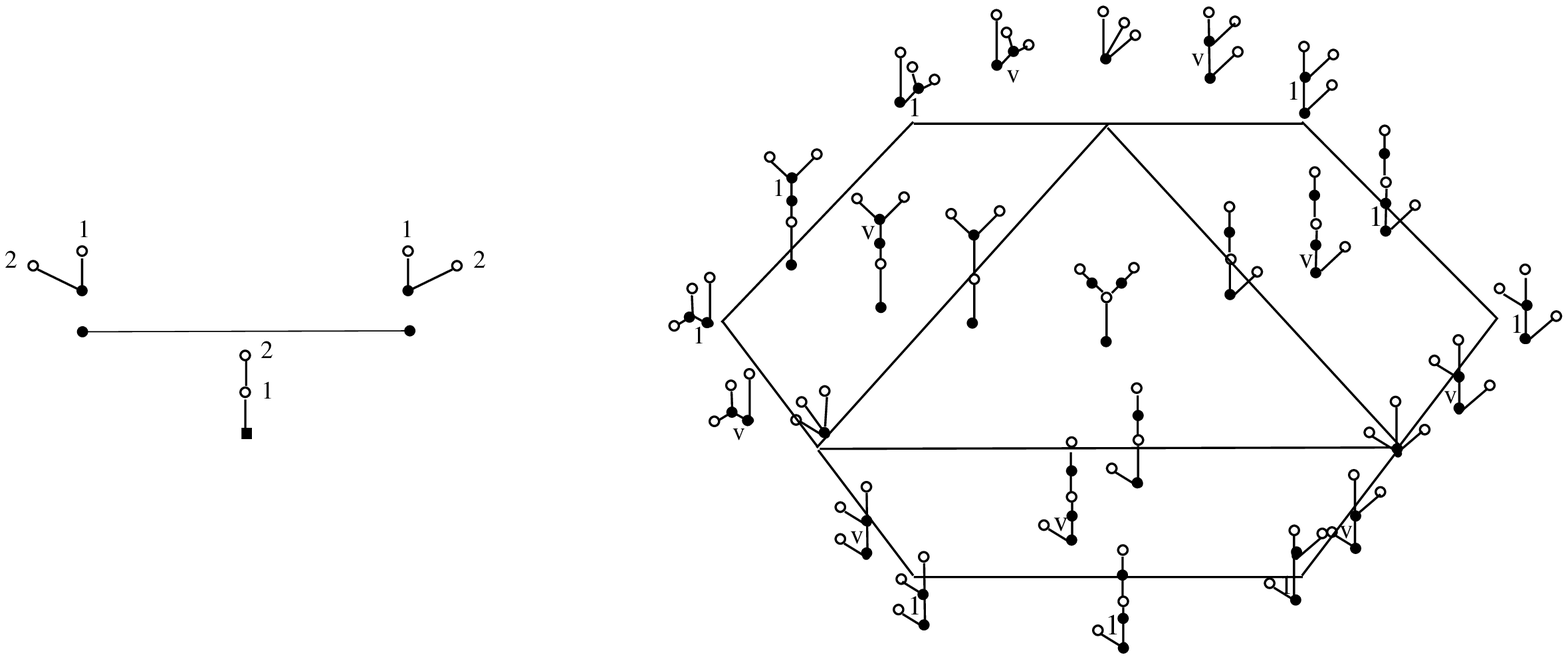}
\caption{\label{W3decomp} The decomposition of $W_2$ and $W_3$.
The labels for the white vertices of $W_3$ are omitted. As depicted,
the special vertex labelled $1$ is always the lowest white
vertex on the center ``stem'' of the tree}
\end{figure}

\subsubsection{A topological realization}
Let $\Cyctop(n)$ be the set of pairs $(\t,\topheight)$ where $\t\in \Tcyc$ is
one of the trees above with $n$ white vertices and
$\topheight:E(\t)\to \Rp$ which satisfy

\begin{enumerate}
\item For all $e\in \Eblack$, $\topheight(e)\leq 1$
\item For all $\a\in \wangles:
 \sum_{\a\in \wangles(v)}\topheight(e)=1$
\end{enumerate}
For convenience, we extend $\topheight$
to all angles at white vertices by marking those
that only have one flag by $1$.
This set obtains a topology induced by collapsing angles and contracting
edges whose $\topheight$ goes to zero.
It is clear that this realizes the cell complex and hence:

\begin{prop}
$\Cyctop$  is a topological PL realization of $W_n$ for
the new CW decomposition and the original CW decomposition.
\end{prop}
\qed

\subsection{Contracting the associahedra and cyclohedra}
\label{flowonesection}

 There is a flow on the two realizations which
contracts all black edges; for $0\leq
t<1:\Psi(t)((\t,\topheight))=(\t,\psi(t)(\topheight))$ where
\begin{eqnarray}
(\psi(t)(\topheight))(\a)&=&\topheight(\a)\text{ for } \a\in \wangles\nn\\
(\psi(t)(\topheight))(e)&=&(1-t)\; \topheight(e)\text{ for } e\in \Evar\nn, 0\leq t< 1
\end{eqnarray}
and $\Psi(1)(\t,\topheight)=(\tilde \t,\topheight|_{\tilde \t})$ where $\tilde
\t$ is the tree $\t$ with all black edges contracted and $\tilde
\topheight$ is $\topheight$ restricted to $\tilde \t$, that is restricted
to the white angles, which remain ``unchanged'' during the construction.
Here ``unchanged'' means that the sets are in natural bijection and
we use this bijection to identify them.

\begin{lem}
The flow contracts $\Cyctop(n)$ to $\Delta^n$ and $K_n$ to a
point and establishes homotopy  equivalences,  actually strong
deformation retracts, between these pairs of spaces.
\end{lem}

\begin{proof}
Using the previous descriptions of the polytopes involved, it is clear that $\Psi$ gives a flow whose image is the purported
one.
\end{proof}

\section{Three CW models, $\CWcact$,$\CWass$ and $\CWheight$, for
the little discs
and their relations }
\label{cellsection}

\subsection{Three CW models}
The basic idea is to form products of the polytopes of the last
section to obtain CW complexes from the various types of trees $\bipartite,
\stable,\heights$.
For $\bipartite$ this has been done in \cite{del}, which
is what we first recall.

\subsubsection{The model $\CWcact$ a.k.a.\ $\Cact^1$}
\begin{df}\cite{del}
We define the CW complex $\CWcact(n)$ to be the following CW complex.
The $k$--cells
are indexed by $\t\in \bipartite(n)$
with $\sum_{v\in \Vwhite(\t)}\arity{v}=k$. The cell corresponding
to a tree is defined to be
\begin{equation}
\label{bipartcelleq}
\cell(\t):=\times_{v\in \Vwhite}\Delta^{\arity{v}}
\end{equation}
The attaching maps are given by using the differential $\del$ on $\bipartite$:
$\del(\cell(\t))=\cell(\del(\t))$ where
we use the orientation and signs dictated by the ordering in equation
(\ref{bipartcelleq}).
\end{df}

\begin{rmk}
This complex was called $\Cact^1(n)$ in \cite{cact,del}
\end{rmk}

The elements in this CW complex are pairs $(\t,\topheight)$ where
$\t\in \bipartite$ and $\topheight$ is a topological ``height'' or
``weight'' function as in \S \ref{topheightsection}; that is a
function $\topheight:\wangles\to (0,1]$ such that $\forall v\in
\Vwhite: \sum_{\a\in \wangles(v)}\topheight(\a)=1$.
 Note that there are no black
edges. The main theorem concerning this complex is:

\begin{thm}\cite{cact,del}
\label{cactthm}
 $|\CWcact|$ is a quasi--operad which
induces an operad structure on $CC_*(\CWcact)$ which in turn is a chain
model for the little discs.
\end{thm}

\subsubsection{The model $\CWass$, a CW realization of $\M$}
\begin{df}\cite{del}
We define the CW complex $\CWass(n)$ to be the following CW complex.
The $k$--cells
are indexed by $\t\in \stable(n)$
with
$\sum_{v\in \Vwhite(\t)}\arity{v}+\sum_{v\in \Vblack} (\arity{v}-1)=k$.
The cell corresponding
to a tree is defined to be
\begin{equation}
  \label{stablecelleq}
\cell(\t):=\times_{v\in \Vwhite} W_{\val{v}} \times \times_{v\in \Vblack}
K_{\arity{v}}
\end{equation}
The attaching maps are given by using the differential $\del$ on $\stable$:
$\del(\cell(\t))=\cell(\del(\t))$ where
we use the orientation and signs dictated by the ordering in equation
(\ref{stablecelleq}).
\end{df}

\begin{lem}
The complexes $\M(n)$ and $CC_*(\CWass(n))$ are isomorphic over $\mathbb Z$.
\end{lem}

\begin{proof}
By construction the two Abelian groups are isomorphic.
Their dg--structures are also compatible by the combinatorics
of the previous section and the construction.
Explicitly, the boundary of  cell is given by
\begin{eqnarray}
\del (\Delta(\t))&=&\sum_{v\in \Vwhite}\pm\del
W_{\val{v}}\times \times_{v'\in\Vwhite\setminus \{v\}}
W_{\val{v'}}\times \times_{v''\in V_{\rm black}} K_{\arity{v''}}\nn\\
&&+ \sum _{v\in \Vblack}\pm \times_{v'\in V_{\rm white}}
W_{\val{v'}}\times \del K_{\arity{v}}\times \times_{v''\in
V_{\rm black}\setminus\{v\}} K_{\arity{v''}}
\end{eqnarray}
where now each summand corresponds to inserting an edge which
is mixed for the first sum and black for the second sum.
This shows that $\M(n)$ and $CC_*(\CWass)(n)$ are
isomorphic complexes.
\end{proof}

\begin{df}
\label{operadassdef}
The {\em induced operad structure} on $CC_*(\CWass):=\{CC_*(\CWass(n))\}$
is the one induced by the isomorphisms $CC_*(\CWass)\cong \M$.
\end{df}

\subsubsection{A new mediating model $\CWheight$}
\begin{df}\cite{del}
We define the CW complex $\CWheight(n)$ to be as follows.
The $k$--cells
are indexed by $(\t,\height)\in \heights(n)$
with
$\sum_{v\in \Vwhite(\t)}\arity{v}+|\Evar|=k$.
The cell corresponding
to a tree is defined to be
\begin{equation}
\label{heightcelleq}
\cell(\t):=\times_{v\in \Vwhite} \Delta^{\arity{v}}
\times I^{\Evar}
\end{equation}
The attaching maps are given by using the differential $\del$ on $\heights$:
$\del(\cell(\t))=\cell(\del(\t))$ where
we use the orientation and signs dictated by the ordering in equation
(\ref{heightcelleq}).
\end{df}

\begin{lem}
\label{topheightslem}
Each element of $|\CWheight(n)|$ corresponds to a pair $(\t,\topheight)$ with
$\t$  a $\{1,\dots, n\}$--labelled stably bipartite tree
and ``heights/weights'' given by
$\topheight:\Eblack(\t)\cup \wangles \to (0,1]$
with the condition that
 $\sum_{\a\in \wangles(v_w)}\topheight(\a)=1$ for all $v_w\in V_{\rm white}$:
\end{lem}

We will call the set of all these pairs $\topheights$.

\begin{proof}
Any element $p$ of $|\CWheight|$ lies inside a unique maximal cell. This
corresponds to a tree $\t\in \heights$. For a black edge $e\in \Evar(\t)$,
we can thus define $\topheight(e)$,
to be the coordinate of $p$ in the factor $I^e$ in $\cell(\t)$,
for the black edges of $\t$ of height $\height(\t)=1$ we set
$\topheight(e)=1$,
and for $\a\in \wangles(v)$,
$\topheight(\a)$ to be given by the barycentric coordinates on
$\Delta^{\arity{v}}\subset \R^{\val{v}}$.
\end{proof}

\subsubsection{Quasi--Operad structure on $|\CWheight|$}
Just as for $|\CWcact|$ above, we can define a quasi--operad
structure
 on the topological level, that is on $|\CWheight|$ which
induces an operad structure on the chain level. We achieve this via an
arc interpretation to realize the space basically
as a sub--quasi--PROP of the
Sullivan--quasi--PROP \cite{hoch1}. We recall
that a topological quasi--operad or quasi--PROP
only has to be associative up to homotopy
(see \cite{cact} for the definition of quasi--operad
and \cite{hoch1} for the definition of quasi--PROP).

\begin{prop}
$|\CWheight|$ is a topological quasi--operad and the quasi--operad maps are
cellular and  induce an
operad structure on $CC_*(\CWheight)\simeq \Z\heights$.
\end{prop}

\begin{proof}
See Appendix \ref{arcapp}, Proposition \ref{appheightprop}.
\end{proof}

\subsection{The relations between the three complexes}
\subsubsection{$\CWheight$ is a refinement of $\CWass$}
\label{treecutsection}
\begin{prop}
\label{refinementprop}
$\CWheight$ is a refinement of $\CWass$, i.e.\  they have
the same realization, and each cell of $\CWheight$ is contained
in a unique cell of $\CWass$.
\end{prop}

\begin{proof}
To show that $|\CWass|\simeq|\CWheight|$ we notice
that each point $p\in |\CWass|$ lies in a
unique maximal cell indexed by a stable tree $\t\in \stable$.
Each cyclohedron $W_{\val{v}}$
or associahedron $K_{\arity{v}}$ appearing as a factor indexed by a
vertex $v$ of $\cell(\t)$ has a decomposition as in
\S\ref{polysection} and our element $p$ lies inside  a unique one of
these finer cells. These finer cells are  of the type $\Delta^k\times I^l$
and are indexed by a tree with heights $\tilde \t(v)\in \heights$,
for each vertex $v$. The coordinates in these cells uniquely
determine the projection to the appropriate factor of $\cell(\t)$
corresponding to the factor $W_{\val{v}}$ or   $K_{\arity{v}}$. To
obtain a pair $(\tilde \t,\topheight)\in \topheights$ as in Lemma
\ref{topheightslem}, we proceed as follows. 
Now for each non--leaf $v$ insert the tree $\tilde\t(v)$
into the vertex $v$. The result is  a stably bipartite tree $\tilde
\t$. We define the function $\topheight$ to is given by the
coordinates of $p$ w.r.t.\ the $\cell(\tilde \t(v))$ for the white
angles and the new black edges  and the markings $1$ for the black
edges stemming from the orignal tree.

Vice versa, given a point $p\in |\CWheight|$, that is a pair
$(\t,\topheight)$, we claim that we can identify it with a point in
one of the finer cells in the decomposition of $\CWass$ above. The
cell of $\CWass$ this point lies in will be indexed by the tree
obtained 
contracting all non--leaf, non--root edges of $\t$ which are not
labelled by $1$ and forgetting the function $\topheight$. Each
pre--image of a vertex, after adding white leaves, will be of the
type $\Tassheight$ or $\Tcycheight$ with a compatible topological
height function $\topheight$. By the previous paragraph this
uniquely determines a point in $|\CWass|$.

It is easy to see that these maps are homeomorphisms that are inverses of
each other. It then follows from the definition of the maps that
 each cell of $\CWheight$ is contained
in a unique cell of $\CWass$.
\end{proof}

\begin{figure}
\epsfxsize = .8\textwidth \epsfbox{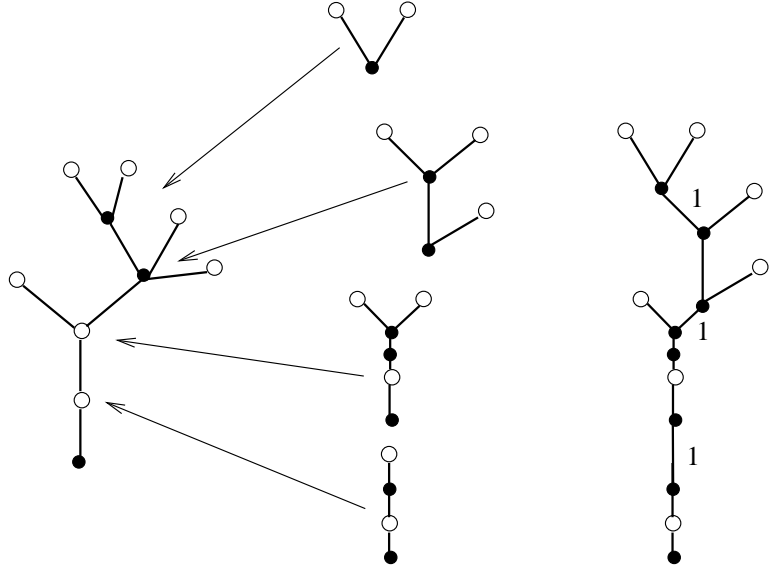}
\caption{\label{treecut} Replacing vertices by trees.}
\end{figure}

For an example of the above procedure see Figure \ref{treecut}.

Using this Proposition, the operad structure on $CC_*(\CWass)$ which
was introduced via the {\it ad hoc} definition \ref{operadassdef}
above can now be induced for the topological level. In other words,
it can be be replaced Proposition \ref{prop1}, which in its precise
form reads:

\begin{prop}
\label{textpropone}
The operad structure of $CC_*(\CWheight)\simeq \Z\heights$
pulls back to $\M\simeq \Z\stable\simeq CC_*(\CWass)$ and
this operad structure coincides with the induced operad structure of
Definition \ref{operadassdef}.
\end{prop}

\begin{proof}
See Appendix \ref{arcapp}.
\end{proof}

\subsubsection{$|\CWheight|$ contracts to  $|\CWcact|$}

Using basically the same flow as in \S\ref{flowonesection}, but now
extended to all of $|\CWheight|$, that is pairs $(\t,\topheight)$, we
can give an explicit deformation retraction.
\begin{df}
We define the flow $\Psi: I \times |\CWheight| \to |\CWheight|$ by
$0\leq t<1:\Psi(t)((\t,\topheight))=(\t,\psi(t)(\topheight))$ where
\begin{eqnarray}
\psi(t)(\topheight)(\a)&=&\topheight(\a)\text{ for } \a\in \wangles\nn\\
\psi(t)(\topheight)(e)&=&(1-t)\; \topheight(e)\text{ for } e\in
\Evar\nn, 0\leq t< 1
\end{eqnarray}
and $\Psi(1)(\t,\topheight)=(\tilde \t,\topheight|_{\tilde \t})$
where $\tilde \t$ is the tree $\t$ with all black edges contracted
and $\tilde \topheight$ is $\topheight$ restricted to $\tilde \t$,
that is restricted to the white angles, which remain ``unchanged''
during the construction. Here ``unchanged'' again means that the
sets are in natural bijection and we use this bijection to identify
them.
\end{df}

\begin{df}
We define $\intop:|\CWcact(n)|\to |\CWheight(n)|$ by mapping a pair
$(\t,\topheight)$ giving a point in $|\CWcact|$ to itself, but now
considered as specifying a point in $|\CWheight|$.
\end{df}
This is possible, since a bipartite tree $\t$ is stably bipartite
and since a bipartite tree has no black edges and hence
 $\Eblack(\t)\cup\wangles(\t)=\wangles(\t)$.

\begin{prop}
\label{retractprop} The topological spaces $|\CWheight(n)|$ and
$|\CWcact(n)|=\Cact^1(n)$ are homotopy equivalent and the homotopy
is given by the explicit flow $\Psi$. This  even yields a strong
deformation retract $r(n)$ of onto the image of
$\intop(|\CWcact|(n))$ and a cellular map.
\end{prop}

\begin{proof}
  It is clear that $\Psi$  is a homotopy and
easy to see that it contracts onto the image of $\intop$, which
remains fixed under the homotopy. This yields the desired statement
\end{proof}

\begin{prop}
\label{textprop4} The sequence of maps  $\stinftop(n):
|\CWass(n)|\stackrel{\sim}{\to}
|\CWheight(n)|\stackrel{r(n)}{\to}|\CWcact(n)|$ induces a
quasi-isomorphism of operads $\stinf:\M\simeq CC_*(\CWass))\to
CC_*(\CWcact)$ on the chain level.
\end{prop}

\begin{proof}
First by Proposition \ref{refinementprop} and Proposition \ref{retractprop}
the composition is cellular and hence indeed induces a map on
the cellular chain level.
We see that any cell of $\stable$ is contracted to a
cell of lower dimension as soon as there is a black vertex whose
valence is greater than $3$, so that these cells are sent to zero. This corresponds
to the fact that $\Psi$ contracts all the associahedra to a point. If
the vertices only have valence $3$ then the black subtrees are
contracted onto the image of $\intop$ which yields a cell of the
same dimension indexed by the tree $\stinf(\t)$. Finally
we know by Lemma \ref{stinflem} that
$\stinf$ is an operadic map. Since $\stinftop$ is a homeomorphism
followed by a strong retraction, the map induced in homology is
an isomorphism.
\end{proof}

\begin{thm}
\label{cellthm}
$\CWass$ is a cell model for the little discs operad whose cells are
indexed by $\stable$.
\end{thm}

\begin{proof}
By  Theorem \ref{cactthm}.
$\CWcact=CC_*(Cact^1)$ is an operadic chain
model for the little discs, hence by
the last proposition we may deduce that $\CWass$ also has
this property.
\end{proof}

\section{The $A_{\infty}$-Deligne conjecture}
\label{ainfdelsection}
In this section we give the solution to the above conjecture
using our results combined with  the action
of the minimal operad $\M$ of \cite{KS}.
We first review the this operation briefly.
Recall that given a tree in $\stable(n)$
there is a natural action on the Hochschild complex
by viewing the tree as a flow
chart. In particular given functions $f_1,\dots, f_n$,  the action of $\t\in
\stable(n)$ is defined as follows: first ``insert'' each of the functions
$f_i$ into the corresponding white vertex $v_i$ and then view the tree as a
flow chart using the operations $\mu_l$ of the $\Ainf$ algebra at
each black vertex of arity $l$ and the brace operation
$f_j\{g_1,\dots,g_k\}$ at each white vertex of arity $k$ to concatenate the functions,
where $f_j$
is the function associated to the vertex $v$ and the $g_i$ are the functions which 
are obtained by following the $k$ flow charts above $v$ corresponding to the
$k$ different branches.

The brace operation 
is defined as \cite{getzler,Kad}
\begin{multline}
\label{bracedef}
  h\{g_1,\dots ,g_n\}(x_1,\dots,x_N) :=\\
\sum_{\footnotesize\begin{tabular}{c}
$1 \leq i_1 \leq \dots \leq i_n \leq |h|:$\\
 $i_j + |g_j|\leq i_{j+1}$\end{tabular}}\pm
h(x_1, \dots, x_{i_1-1},g_1(x_{i_1}, \dots, x_{i_1+|g_1|}), \dots, \\
\dots ,x_{i_n-1}, g_n(x_{i_n}, \dots, x_{i_n+|g_n|}), \dots , x_N)
\end{multline}

\begin{thm}[Main Theorem]
There is a cell model $\CWass$ for the little discs operad, whose operad
of cellular chains $CC_*(\CWass)$  acts on
the Hochschild cochains of an $\Ainf$ algebra inducing
the standard operations of its homology on the cohomology.
Moreover, this is minimal in the sense
that  the cells correspond exactly to the natural operations obtained
by concatenating the functions and using the $\Ainf$ structure maps.
\end{thm}

\begin{proof}
This follows from Theorem \ref{cellthm} in conjunction with the Theorem
of \cite{KS} that the operad $\M\simeq\stable$ acts in a dg--fashion on
Hochschild cochains of an $\Ainf$  algebra.
\end{proof}

\renewcommand{\theequation}{\Alph{thesection}-\arabic{equation}}

\renewcommand{\thesection}{\Alph{section}}
\setcounter{equation}{0}  
\setcounter{section}{0}

\app{Connection to arcs and polygons with diagonals}

 In this
Appendix, we give the connection of the CW complexes
 to the arc operad of \cite{KLP} and the
Sullivan--quasi--PROP of \cite{hoch1}.
All of the (quasi--) operad structures we are concerned
with are based on the two mentioned structures, 
and we use these facts to give proofs of Theorem \ref{theorem2}
and Proposition \ref{prop1}.
There are actually three different pictorial realizations
for the same objects: arc graphs, ribbon graphs and
trees.
These correspondences have been worked out in full detail in
\cite{del,hoch1,manfest}, and we will
content ourselves with a brief
review of the salient
features referring the fastidious reader to these papers.

\subsection{The arc picture}
First we would like to recall that an element of $\DArc$ is
 a surface $F_{g,n+1}^r$ of genus
$g$ with $n+1$ boundary components labelled by $\{0,\dots,n\}$
and $r$ punctures
 with marked boundary, that is one marked
point per boundary component together with two sets
of data, an arc graph and weights.

An arc graph is a collection of
arcs, that is embedded curves from boundary to boundary that
\begin{enumerate}
\item Do not hit the marked points.
\item Do not intersect.
\item Are not parallel. This means that they are not homotopic
to each other, where the endpoints  may not cross endpoints of
other arcs or the marked points.
\item Are not parallel to the a part of the boundary, where these now include the marked
points.
\item All boundaries are hit, that is they have at least one incident arc.
\end{enumerate}
considered up to the action of the pure mapping class group that
keeps all punctures and marked points pointwise fixed
and the boundaries setwise fixed.

Weights for an arc graph are given by assigning a weight to
each arc, that is a map from the set of all arcs to $\Rp$.
We will only need to consider $g=r=0$ in the present discussion
and we {\em fix this from now on}.

\subsubsection{Gluing in the arc picture}
The gluing is understood as a gluing of partially measured
foliations, which can be paraphrased as follows. Realize the
arcs with weights as bands with width. If two sets of bands incident
to two boundaries have the same total width, just splice them together along
their leaves. That is glue the bands and cut along the common partition.

The different operad/quasi--operad/quasi--PROP structures \cite{KLP,hoch1}
are basically built in the same fashion. First pick two boundaries
to be glued, then scale such that the weights agree, and finally
glue the boundaries and the foliations as explained above.
We will have a new feature for $|\CWheight|$ since
the topological gluing will involve a forth step
of renormalizing.

Regardless of this there are two parts to the gluing, one combinatorial,
where the combinatorics govern the types of arcs that occur and the
second topological, which is the part dictated by the particular weights.
On the chain level, we only want to keep the combinatorics.

\subsection{Embedding $|\CWheight|$ into $\DArc$ and Generalized Cacti}
\label{linfsection}
Just as there is a topological embedding of $\Cact^1$ into the arc
operad $\Arc$ of \cite{KLP}, there is also such an embedding of $|\CWheight|$
into $\Arc$. We let $\Linf$ be the subspace which consists of
those arc families that  are of genus 0 with no punctures, arcs 
running only from $i$ to $0$ {\em and possibly} arcs running from $0$ to $0$,
which satisfy the following conditions. There is a representative of
projective weights on the arcs such that
\begin{enumerate}
\item No arc running from $0$ to $0$ homotopic to a boundary $i$
together with one arc from $i$ to $0$ where the marked point
is considered to be part of the boundary.
\item The linear orders at the boundaries $i$ are (anti)--compatible
with the linear order at $0$. That is, if for two arcs $a$ and $b$
which hit the boundary $i$  $a<_0b$ in the order at $0$, then we have $a>_ib$ in the order
at the boundary $i$.
\end{enumerate}

The space $|\CWheight|$ corresponds to
 the subset $\Linf^1\subset \Linf$, which additionally
satisfies
\begin{enumerate}
\setcounter{enumi}{2}
\item The weight of each arc from $0$ to $0$ is $\leq 1$
\item The sum of the weights for each boundary except $0$ is one.
\end{enumerate}

\begin{figure}
\epsfxsize = \textwidth \epsfbox{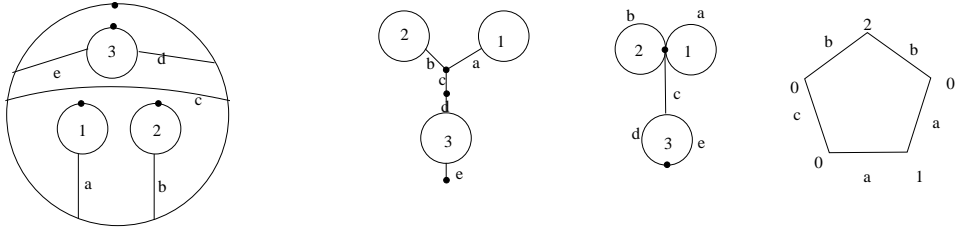}
\caption{\label{simpleex} An arc graph, its tree,  cactus
representation and one of its polygons}
\end{figure}

In the following, we give a brief  translation primer for
the different combinatorial pictures. An example
 is given in Figure \ref{simpleex}.

\subsubsection{From Arc graphs to ribbon graphs}
Given an arc family in $\Arc$ we first define its dual ribbon
graph. This has one vertex for each complementary region and one
edge for between the two (not necessarily distinct) regions on the
different sides of each arc. See \cite{cact,del} for more details.
Every cycle of the ribbon graph corresponds to exactly one boundary component.
Since the boundary components were oriented and marked, the ribbon graph will be marked as well, that is, there is one distinguished flag in each cycle
that points in the direction of the orientation and has its vertex in
the region that contains the marked point.

Notice that in our case, since all arcs run to zero, there is a distinguished
cycle which
runs through all the edges. That is, the ribbon graph is tree-like
in the terminology of \cite{hoch1}. In this correspondence
each arc corresponds to an edge, and hence if the arcs have weights, so have
the edges.

\subsubsection{From ribbon graphs to trees}
For a tree--like ribbon graph, define its
{\em incidence graph}   to be given
by
one white vertex for each cycle excluding the distinguished one and
a black vertex for each previous vertex, where we join two black vertices
if they are joined in the original graph along an edge which does
not belong to the non-distinguished cycles and we join a white and a
black vertex if the black vertex lies on the cycle given by the
white vertex. The tree is rooted and planted by taking the flag
corresponding to the  marked
flag of the graph as the marked flag of the tree.
Now the edges correspond to the white angles and the black edges and
hence these carry the weights.

\subsubsection{From $\heights$ trees to ribbon graphs}

Given a tree in $\heights$ we first ``blow--up'' the white
vertices to cycles and then contract all the images of the mixed edges.
In the blowing up process each angle becomes an edge of the ribbon graph with
the two flags of the angle incident to the two vertices of the new edge preserving their orders.
The labels are now on all of the edges.

\subsubsection{From ribbon graphs to arc graphs}
It is well known that thickening a ribbon graph gives rise to a surface
with an embedding of the ribbon graph as the
spine. Taking the dual graph on
the surface basically yields an arc graph. For the missing 
makings, we 
mark the respective boundary of the respective
region containing the marked flag of the cycle.
The weights pass along the bijection of the edges and
the markings.  We refer to
\cite{del} for more details.

\subsubsection{Description of $\heights$ in terms of polygons}
By the above procedure every tree in $\heights$ translates
to an element in $\DArc$. Cutting along the arcs
decomposes the surface into pieces, and, as we
fixed that $g=s=0$ above, these pieces are polygons. These polygons are $2n$-gons with sides
alternatingly corresponding to pieces of the boundary and arcs.
We obtain a set of polygons by
contracting all sides corresponding to boundaries
and call these the complementary polygons.

We have the following translation table

\noindent\begin{tabular}{|@{\extracolsep{0pt}}l|l|}
\hline
$\heights$&$\DArc$\\
 \hline\hline
mixed edge&arc from $0$ to $i\neq 0$\\
\hline
black edge&arc from $0$ to $0$\\
\hline
 There are no white edges& the tree is an intersection graph\\
\hline
 There are no black vertices of valence $2$ &no parallel arcs\\
both of whose edges are black.&\\
There are no black vertices of valence $2$&there are no triangles among the
\\
 with one edge black and
the other edge& complementary polygons, where two \\
 a leaf edge {\em unless} the vertex is the root.& edges correspond to the same arc.\\
\hline
Trees obtained by cutting &complementary regions of the\\
black edges marked by $1$&arcs from $0$ to $0$ of weight $1$.\\
\hline
\end{tabular}

\subsubsection{Generalized Spineless Cacti}
Yet another way to picture the trees is to look at the ribbon graph
as a new version of cacti. Here one is now allowed to have edges between
the lobes. We  define  $\Cinf$  to be the
space of metric marked ribbon graphs corresponding
to the subspace $\Linf$ of $\DArc$.
\begin{prop}
$\Linf$ is a sub-operad and hence $\Cinf$ is an operad.
\end{prop}
\begin{proof}
The claim boils down to checking that the conditions of $\Linf$ are
stable under gluing, which they are.
\end{proof}

We also let $\Cinf^1$ be the  space of ribbon graphs corresponding
to $\Linf^1$.
\subsubsection{Gluing in $\Cinf$}
The gluing operation defined above is reminiscent of the definition of the gluing
of $\Cact$ as defined in \cite{cact}.
If we are given two generalized normalized spineless
 cacti $c_1,c_2\in \Cinf^1$ then $c_1\circ_i c_2$
is the generalized normalized spineless cactus obtained as follows.
Glue $c_2$ into the cycle $i$ of $c_1$ by identifying the cycle $0$ of $c_2$
with the cycle $i$ of $c_1$, where these cycles are considered to be
parameterized over $S^1$ by the metric on their edges and their marked points.
Here it is important, that we scale the total length, i.e.\ the sum of
weights of all the edges, of  $\t'$  to fit the
the sum of the weights of the edges of the lobe $i$ of $\t$. For the quasi--PROP
structure, we will scale the other way around, that is scale the lobe
to fit. Also to fit the combinatorics, we will need to renormalize
this construction.

\subsection{The Sullivan--quasi--PROP of \cite{hoch1}}
We briefly review the  Sullivan--quasi--PROP of \cite{hoch1}, but
refer the reader to {\it loc.\ cit.} for details.

In order to make contact with the quasi--PROP structure, we
need to additionally assume that the boundary labels of
the surfaces in question are divided into $\In$ and $\Out$ boundaries
with labels. Correspondingly we will obtain spaces $\DArc(\In,\Out)$. If $|\In|=n$ and $|\Out|=m$
this is naturally a collection of $\;\Sn\times \Sm$ modules. We will
simplify and fix $\In=\{1,\dots, n\}$, $\Out=\{1,\dots,m\}$.

We let $\Dsul$ be the collection of subspaces of the spaces of
$\DArc(\In,\Out)$ in which there are only arcs running from the $\In$ to
the $\Out$ and possibly
from the $\Out$ to the $\Out$
boundaries  {\em and} there is no empty $\In$ boundary.
This space was denoted $\overline{\DArc}^{i\nleftrightarrow i}$
in \cite{hoch1}.
We define $\CWsul\subset \Dsul$ to be the subspace of graphs
whose sum of weights of arcs incident to every $\In$ boundary
vertex is one and whose arcs
from $\Out$ to $\Out$ have weights $\leq 1$. This is naturally a CW complex.

In \cite{hoch1} we defined the quasi--PROP compositions
on $\Dsul$
by scaling
the input $i$ individually to the weight of the output $j$ it is
glued to. This yields topological quasi--PROP structure $\bullet_{i,j}$.
 Notice that in the gluings
 one only scales at the $\In$ boundaries which are
to be glued so that the
weights on the $\In$ boundaries which remain after gluing are unchanged
as are the weights of the arcs from $\Out$ to $\Out$ boundaries.
Hence $\CWsul$ is a sub--quasi--PROP.

\begin{prop}
\label{hochprop}
The compositions $\bullet$ define a homotopy--PROP structure on
the cell complex $\CWsul$.
\end{prop}
Here homotopy--PROP means a PROP that is associative up to homotopy
\cite{hoch1}.

\begin{proof}
First the fact that $\CWsul$ is a cell complex follows in the previous
pattern. The cells are just indexed by the relevant graphs.
It is clear that $\CWsul$ is stable under composition
\end{proof}

Although the PROP structure $\CWsul$  is cellular,
 it does not directly
yield exactly the $dg$--PROP structure we are looking for.
To make the proofs simpler we again restrict to $g=s=0$ and
deal only with the special sub--structure we are interested in.
Namely, we consider $\Linf(n)$ as a subspace of $\CWsul(n,1)$ if we declare $0$
to be in $\Out$ and all other inputs to be in $\In$. We will identify $\Cinf$
with $\Cinf$ and we will also use the term lobe for a cycle corresponding
to an $\In$ boundary.

We will also call an arc black if it runs from $0$ to $0$ as it
will give rise to a black edge and we will call the other arcs white arcs,
as they will give rise to white angles.

\subsubsection{Renormalized Gluing in $\Linf^1$}
In the gluing procedure of the quasi--PROP given by $\bullet$,
black bands might be split and although this will induce the right
kind of combinatorics on the topological level, it actually yields
the wrong type of combinatorics on the chain level. This is simply due to the fact
that after splitting a band it can never have weight $1$.
In order to rectify the situation, we define a slightly modified gluing procedure $\bar\bullet$
as follows. First glue according to $\bullet$ and then for each black arc that
is split into $n$ arcs we rescale according 
to the radial projection
 $\Delta^{n-1}\to I^n$ that maps the simplex homeomorphically   
to the faces
of $I^n$ which have one or more entries equal to $1$.
 To be precise, if the black arc that is split has weight $w$
and the $n$ arcs it splits into have weights $t_1,\dots t_n$ with
$\sum t_i=W$ then we re-scale the weights
to $(\tilde t_1,\dots \tilde t_n)$,
which is the the image of
$(t_1,\dots, t_n)$ under the radial projection onto the cube $[0,W]^n$.

\begin{prop}
\label{appheightprop}
$\Linf^1$ is a sub--CW complex of $\CWsul$ and hence a CW complex.
The operations $\bar \bullet$
endow $\Linf^1$ with a topological quasi--operad structure,
which is equivalent as a quasi--operad to its
topological sub--quasi--PROP structure.

Furthermore, the operations $\bar \bullet$ induce an operad structure
on  $CC_*(\Linf^1)$ and moreover $CC_*(\Linf)\simeq \Z\heights$.
The same statements hold true for $\Cinf^1$, by identifying it with $\Linf^1$
\end{prop}

\begin{proof}
It is clear that $\Linf^1$ is a sub--CW complex and stable
under the quasi--PROP compositions. The difference between $\bullet$
and $\bar \bullet$ is the radial projection which is homotopic to the
identity and hence the two structures are both associative up to homotopy
and this homotopy gives the equivalence.

Now by taking the intersection graph of a ribbon graph, we see that additively
$CC_*(\Linf^1)=CC_*(\Cinf^1)=\Z\heights$. Taking
the composition $\bar \bullet$
means that indeed we are allowed all the combinations of putting
branches into the angles and into the black edges. The former corresponds
to the splitting of a white arc and the latter to the splitting of
a black arc. Now $\bar \bullet$ was chosen so that inserting
into a black edge gives exactly the summands corresponding to the distribution
of labels.  It is now straightforward to check that these gluings
are now strictly associative.
\end{proof}

\subsubsection{Sub--Quasi--PROP Structure of $|\CWass|$ and $|\CWheight|$}
\begin{thm}[Theorem \ref{theorem2}]
\label{texttheorem2}
The realizations $|\CWass|\simeq |\CWheight|$ and $|\CWcact|$ are
all topological quasi--operads and sub--quasi--PROPs of
the Sullivan--PROP $\CWsul$.
There is also a renormalized quasi--operad structure such that
the induced quasi--operad structures on their cellular chains
$CC_*(\CWass)\simeq\Z\stable$, $CC_*(\CWheight) \simeq \Z\heights$
and $CC_*(\CWcact)\simeq\Z\bipartite$ are operad structures
and coincide with the respective combinatorial operad structure on the trees.
Moreover, all these operad structures are models for the little discs operad.
\end{thm}

\begin{proof}[Proof of Theorem \ref{texttheorem2} and
Proposition \ref{textpropone}]
Taking the intersection graph of the elements of $\Cinf^1$ we obtain precisely
$|\CWheight|$ so that the claims for $\CWheight$ follow from
Proposition \ref{hochprop} and \ref{appheightprop}. Now by the cellular map that identifies $|\CWass|$
with $\CWheight$, each cell of $\CWass$ is a sum of cells of $\CWheight$.
What we must show is that composing sums of these cells indeed gives a
sum of cells. This is most easily demonstrated using $\Cinf^1$.
In this language the argument is analogous
to the one in \cite{del}. Explicitly we claim that if $c_1$ and $c_2$ are
elements of a fixed cells $\cell(\t_1)$ and
$\cell(\t_2)$ of $\CWass$, that is, a sum of cells of $\CWheight$,
as they vary throughout these cells $c_1\circ_i c_2$ produces
exactly  the elements of the cells corresponding to the tree
$\t_1\circ_i\t_2$. This is obvious if one considers $c_2$
as a subgraph of $c_1\circ_i c_2$ whose white
vertices are re-labelled according to the operad composition.
This then allows to extract
$c_1$ and $c_2$ from the data  and $c_1\circ_i c_2$ uniquely after we
fix the number of lobes of $c_1$ and $c_2$ and include these and $i$ into
the data as well.
Hence looking at a possible configuration in $\cell(\t_1\circ_i\t_2)$
we see that it comes precisely from one $c_1$ and $c_2$ via $\circ_i$.
This proves the claims about the chain level of $\CWass$ in
Theorem \ref{texttheorem2} and \ref{textpropone}.

On homology all these models induce the same structure. The map
$\stinf$ is operadic and the same is true for the one
induced by the retraction. On homology the operad structure  is known
by \cite{cact} to be isomorphic to the homology of the little discs operad.
\end{proof}

We can actually also prove a little more:
\begin{thm}
\label{extrathm}
$|\CWass|\simeq |\CWheight|$ are equivalent as topological
quasi--operads to the sub-operad $\Linf$ which in turn is equivalent to the
little discs operad.
\end{thm}

\begin{proof}
It is clear that $\Linf$ is a sub-operad of $\DArc$. For
both $\Linf$ and $\Linf^1$, we can simultaneously
scale to length $0$ all the edges running form $0$ to $0$.
This gives a homotopy equivalence of $\Linf$ with the model $\Cact$
for the little discs operad (see \cite{cact}) and of $\Linf^1$ with the equivalent
model $\Cact^1$.
Furthermore, if for $\Linf$
we also scale the weights on the other edges at the same time, so that they sum up
to $1$ on each boundary we can directly contract it to $\Cact^1$. 
Another way to see the homotopy equivalence of 
$\Linf$ and $\Linf^1$ is to notice that the sum of the weights on the boundaries $1,\dots,n$
contributes
a contractible factor of $\Rp^n$. Hence we have homotopy equivalences of both spaces with $\Cact^1$
and it is a straightforward check that this is through homotopies of quasi--operads.
This can be done analogously to the
argument for $\Cact^1$ relative to $\Cact$ given in \cite{cact}.
Hence both are equivalent to $\Cact^1$ and thus to each other and the little
discs operad (as quasi--operads).
\end{proof}

\app{Sequential Blow--ups/downs for the Cyclohedron}
The subdivision of the cyclohedron by the trees with height $\Tcycheight$
give us an explicit way to blow up the simplex. For this
we notice that the number of black edges marked by $\var$ gives
a depth function $\depth(\t)=|\Evar|$. In the top--dimensional
cells of $W_n$ $\depth(\t)+\val{v}=n$. Here $v$ is the special vertex labelled
by $1$ that is allowed to be a non--leaf vertex.

\begin{thm}[Theorem \ref{theorem3}]
\label{texttheorem3}
There is a new decomposition of the cyclohedron $W_{n+1}$ into a simplex and cubes.
Correspondingly, there is an iterated ``blow--up'' of the simplex to a
cyclohedron, with $n-1$ steps. At
 each stage the objects that are glued on are a
 product of a simplex $\Delta^{n-k}$ and a cube $I^k$
where the factors $\Delta^{n-k}$ is attach to
the codimension $k$--faces of the original simplex.
\end{thm}

\begin{proof}
We  use the depth function to index the iteration. There is only one element of
depth $0$ and this corresponds to the simplex. This is step 0 and the starting point
of the iteration.
All trees of higher depth have
a product of a simplex and a cube as their cell. Furthermore, we notice that
for a new edge in $\Evar$ to appear in a tree indexing an adjacent maximal cell, 
we first have to collapse one effective white angle. Hence
we obtain an iteration for the gluing
of the maximal cells, by first collapsing one angle, then allowing
to collapse 2 angles and so on. This iteration according to the number
of angles collapsed is precisely by depth.
Finally, the $\Delta^{n-k}$ factors are naturally identified with the
codimension  $k$ faces of $\Delta^n$ as
they correspond to collapsing $k$ angles and the choices
for these angles are precisely indexed by the different faces;
see \S\ref{simplexsection}.
\end{proof}

In the first step one ``fattens'' the faces of the simplex $\Delta^n$
by gluing a $\Delta^{n-1}\times I$ onto each face and in 
the last step one simply glues in cubes.

We illustrate this for $W_3$ and $W_4$.
The figure for $W_3$ is Figure \ref{W3decomp}, where there is only one blow--up.
\begin{itemize}
\item[{\it Depth 0}.] This is the simplex $\Delta^2$.
\item[{\it Depth 1}.] The new elements are
products $\Delta^1\times I=I^2$.
There are exactly 3 of these which are glued onto the sides of $\Delta^2$.
\end{itemize}
This gives an nonegon, but identifying 3 pairs of sequential sides
and all the top--dimensional cells, we are left with the usual hexagon picture;
see Figure \ref{W3decomp}.

\begin{figure}
\epsfxsize = .8\textwidth \epsfbox{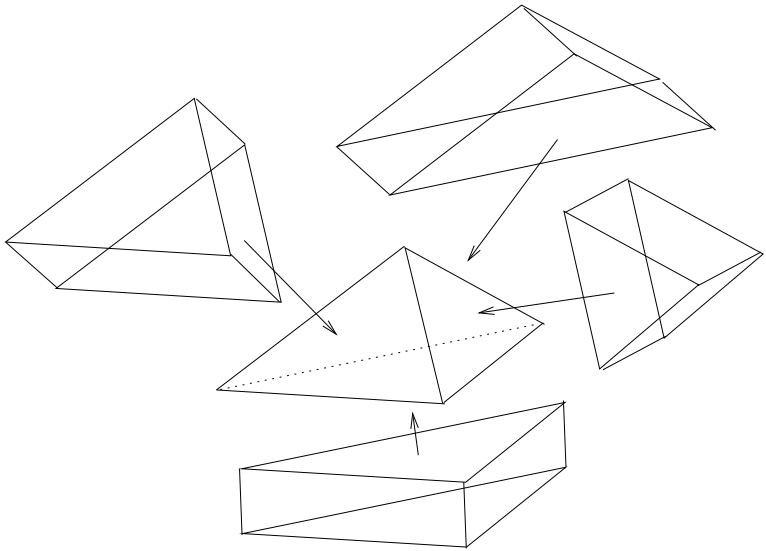}
\caption{\label{W4step1} Step 1: Gluing on four $\Delta^2\times I$s to $\Delta^3$}
\end{figure}

\begin{figure}
\epsfxsize = .8\textwidth \epsfbox{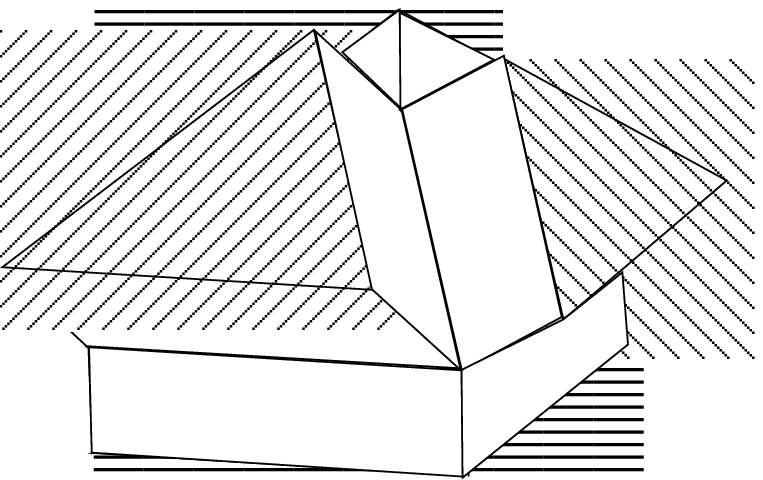}
\caption{\label{W4step2} The result
of the first gluing}
\end{figure}

\begin{figure}
\epsfxsize = .8\textwidth \epsfbox{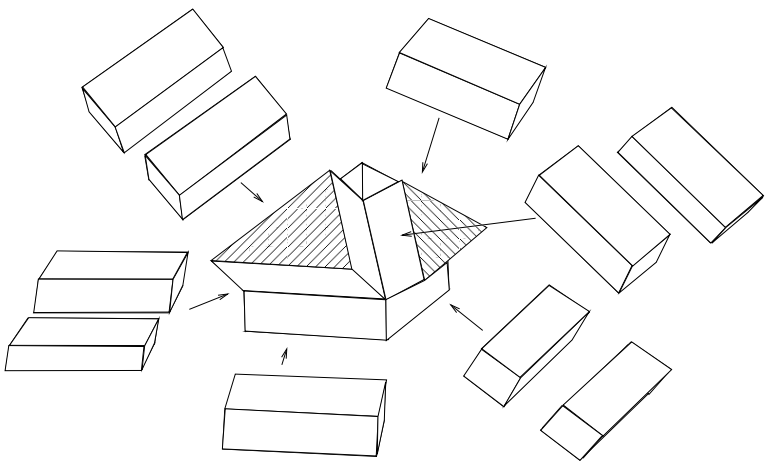}
\caption{\label{W4step2b} Step 2: gluing on 10 $\Delta^1 \times I\times I=I^3$s}
\end{figure}

\begin{figure}
\epsfxsize = .8\textwidth \epsfbox{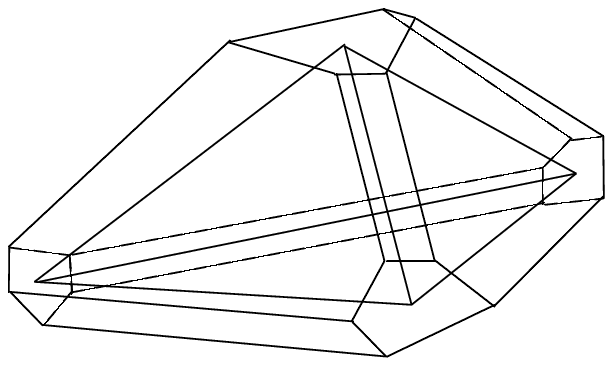}
\caption{\label{del4inW4} The simplex $\Delta^3$
inside $W_4$ after the construction}
\end{figure}

For $W_4$ there are 2 blow--ups and the details are illustrated in Figure \ref{del4inW4}.
\begin{itemize}
\item[{\it Depth 0}.] This is the simplex $\Delta^3$.
\item[{\it Depth 1}.] The new elements are products $\Delta^2\times I$.
There are exactly 4 of these which are glued onto the 4 faces of $\Delta^3$,
see Figure \ref{W4step1}. The result is given in Figure \ref{W4step2}.
\item[{\it Depth 2}.] There are 10 elements of the form $\Delta^1\times I^2=I^3$
which are glued in. This is asymmetric (as it should be). Four of
the edges are associated to two cubes and two of the edges to only one cube.
The latter two edges do not intersect; see Figure \ref{W4step2b}.
\end{itemize}

\begin{figure}
\epsfxsize = .8\textwidth \epsfbox{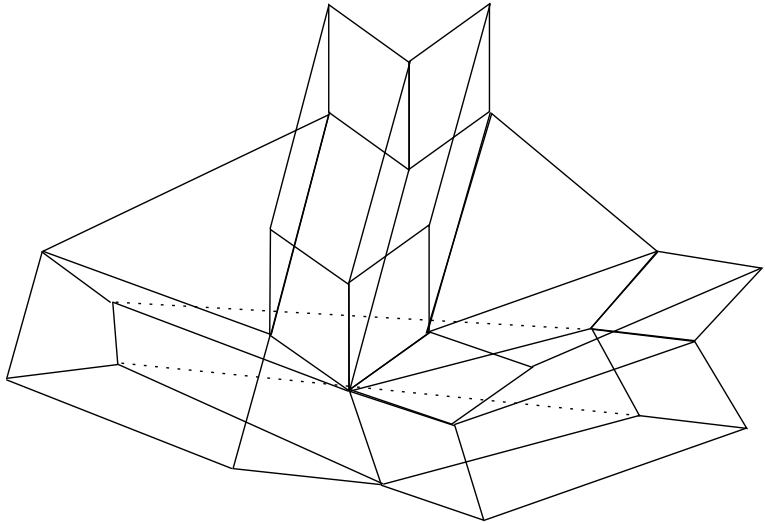}
\caption{\label{starW4} A vertex after the blow--up}
\end{figure}

After the second blow--up, we see that at each vertex there are precisely 2+2+1 cubes,
which effectively replace the vertex  by 5 squares which assemble to a $K_4$;
see Figure \ref{starW4}. If we straighten out the polytope and
consolidate the cells, we obtain the usual picture of $W_4$
(see Figure \ref{del4inW4}) where we now is $\Delta^3$ realized inside
$W_4$.

\begin{rmk}
Notice that the procedure above actually gives a PL embedding of
$W_{n}$ into $\R^{n-1}$.
\end{rmk}

\begin{rmk}
This iteration can also be understood purely in terms of bracketings instead of trees.
We refer the interested reader to \cite{Rachelthesis}.
\end{rmk}

\begin{rmk}
We can alternatively think of the gluings  as a blow--up that
comes about by cutting edges to blow up the faces. In the first step,
we cut along all the edges and then in the second step,
we cut along the four non--special edges. For the purposed of the present
paper it was important however, that we have an explicit embedding of the simplex and
a retraction to it.
\end{rmk}

\end{document}